\documentstyle[amscd]{amsart}
\newsymbol\restriction 1316
\newsymbol\varkappa 207B
\def\dj{d\kern-.30em\raise1.25ex\vbox{\hrule width .3em height .03em}}
\def\Dj{D\kern-.75em\raise0.75ex\vbox{\hrule width .3em height .03em}
\kern.03em}
\makeatletter
\renewcommand{\subsection}{\@startsection{subsection}{2}{\z@}%
{\baselineskip}{0.5\baselineskip}{\defaultfont\bf}}
\makeatother
\def\im{\mathrm{im}}
\def\bla#1{$(${\it #1\/{}}$)$}
\def\e{\epsilon}
\def\k{\kappa}
\def\K{\varkappa}
\def\id{\mathrm{id}}
\def\ad{\mathrm{ad}}
\def\inv{i\!\hspace{0.8pt}n\!\hspace{0.6pt}v}
\def\restr{{\restriction}}
\def\rtw#1{#1^{\mathrm{r}}}
\def\ltw#1{#1^{\mathrm{l}}}
\def\ilG{\iota_\Gamma^{\mathrm{l}}}
\def\Ginvr{{^{\inv}\Gamma}}
\def\sinvr{{_*\sigma}}
\def\cirr{\bullet}
\def\pir{\varsigma}
\def\adj{\varpi}
\def\lac{\ell_\Gamma}

\def\rac{\varrho_\Gamma}
\def\irG{\iota_\Gamma^{\mathrm{r}}}
\def\mlG{m_\Gamma^{\mathrm{l}}}
\def\mrG{m_\Gamma^{\mathrm{r}}}
\newtheorem{lem}{Lemma}[section]
\newtheorem{pro}[lem]{Proposition}
\theoremstyle{definition}
\newtheorem{defn}{Definition}
\numberwithin{equation}{section}
\begin{document}
\title[Braided Differential Calculi]{First-Order Differential Calculi
Over\\Multi-Braided Quantum Groups}
\author{Mi\'co {\Dj}ur{\Dj}evi\'c}
\address{Instituto de Matematicas, UNAM, 
Area de la Investigacion Cientifica, 
Circuito Exterior, Ciudad Universitaria,
M\'exico DF, CP 04510, MEXICO\newline
\indent {\it Written In}\newline
\indent Faculty of Physics, University of Belgrade, pBox 550,
Studentski Trg 12, 11001 Beograd, SERBIA}
\maketitle

\begin{abstract} A differential calculus  of  the  first  order  over
multi-braided quantum groups is developed. In analogy  with 
the standard theory, left/right-covariant and bicovariant differential
structures are introduced and investigated. Furthermore,
antipodally covariant calculi are studied. The concept of the
*-structure on a multi-braided quantum group is formulated, and in particular
the structure of left-covariant *-covariant calculi is analyzed. 
A special attention is given to differential calculi covariant with respect
to the action of the associated braid system. 
In particular it is shown that the left/right braided-covariance
appears as a consequence of the left/right-covariance relative to the
group action. Braided counterparts of all basic results of the standard
theory are found. 
\end{abstract}

\bigskip
\tableofcontents

\filbreak
\renewcommand{\thepage}{}
\section{Introduction}

The basic theme of this study is the analysis of the first-order differential
structures over multi-braided quantum groups. Standard braided quantum
groups \cite{Maj} are included as a special case into the theory of
multi-braided quantum groups \cite{D-BQG}. The difference between two types
of braided quantum groups is in the behavior of
the coproduct map. In the standard theory, the coproduct $\phi\colon
\cal{A}\rightarrow\cal{A}\otimes\cal{A}$ is interpretable as a morphism
in a braided category generated by the basic algebra $\cal{A}$ and
the associated braiding $\sigma\colon\cal{A}\otimes\cal{A}\rightarrow
\cal{A}\otimes\cal{A}$. In our generalized framework two standard
pentagonal diagrams expressing compatibility between $\phi$ and $\sigma$ are
replaced by a single more general octagonal diagram.

It then turns out that the lack of the functoriality of the coproduct map
is `measurable' by a second braid operator $\tau\colon\cal{A}\otimes\cal{A}
\rightarrow\cal{A}\otimes\cal{A}$. Furthermore, two braid operators
generate in a natural manner a generally infinite system of braid operators
$\sigma_n\colon\cal{A}\otimes\cal{A}\rightarrow\cal{A}\otimes\cal{A}$, where
$n\in\Bbb{Z}$, 
which elegantly express twisting properties of all the maps appearing
in the game. This explains our attribute {\it multi}-braided,
for the structures we are dealing with. 

Multi-braided quantum groups include various completely `pointless'
structures, overcoming in such a way an inherent geometrical inhomogeneity
of standard quantum groups and braided quantum groups. This inhomogeneity
is explicitly visible in geometrical situations in which `diffeomorphisms'
of quantum spaces appear. For example, in the theory of locally trivial
quantum principal bundles over classical smooth manifolds \cite{D-QPB1}
a natural correspondence between quantum $G$-bundles (where $G$ is
a standard compact quantum group) and ordinary $G_{cl}$-bundles
(over the same manifold) holds. Here $G_{cl}$ is the classical part of $G$. 

The multi-braided formalism reduces to the standard braided quantum
groups iff $\sigma=\tau$, which means that all the operators
$\sigma_n$ coincide with $\sigma$. This is also equivalent to the
multiplicativity of the counit map. 

In the formalization of the concept of a first-order calculus,
we shall follow \cite{W}: If the algebra $\cal{A}$ represents a
quantum space $X$, then every first-order calculus over
$X$ will be represented by an $\cal{A}$-bimodule
$\Gamma$, playing  the
role of the $1$-forms on $X$, together with a standard derivation
$d:\cal{A}\rightarrow\Gamma$,
playing the role of the  differential. 
Such a formalization reflects noncommutative-geometric \cite{C}
philosophy, according to which the concept of a  differential form 
should be the starting  point  for  a  foundation  of  the quantum 
differential calculus.

    The paper is organized as  follows. In Section~2 we first study
differential calculi over a quantum space $X$,
compatible in the appropriate sense with a single braid operator
$\sigma\colon\cal{A}\otimes\cal{A}\rightarrow\cal{A}\otimes\cal{A}$. 
In this context, left/right, and bi-$\sigma$-covariant
differential structures are distinguished. The  notion  of  left
$\sigma$-covariance
requires a natural  extendability  of $\sigma$ to  a  flip-over  operator 
$\ltw{\sigma}\colon\Gamma\otimes \cal{A}\rightarrow\cal{A}\otimes \Gamma$.
Similarly, right $\sigma$-covariance  requires extendability 
of $\sigma$ to a flip-over operator 
$\rtw{\sigma}\colon\cal{A}\otimes \Gamma\rightarrow\Gamma\otimes \cal{A}$.
Finally, the concept
of bi-$\sigma$-covariance is simply a symbiosis of the previous two. 

\renewcommand{\thepage}{\arabic{page}}
We shall then briefly analyze general situations in which the calculus is
covariant  relative  to  a  given  braid  system $\cal{T}$
operating in $\cal{A}$.

At the end of Section~2 we begin the  study  of  differential 
structures over multi-braided quantum groups.  We  shall  prove
that if $\cal{A}$ and $\sigma$ are associated to a multi-braided
quantum group $G$ then  the  left/right
$\sigma$-covariance implies the left/right  $\sigma_n$-covariance,  for 
each $n\in\Bbb{Z}$. We shall also analyze interrelations between
all possible flip-over operators and maps determining the group structure.

All considerations with braid operators can be
performed at the language of braid and tangle diagrams, as in the framework
of braided categories \cite{Maj}. The unique additional moment is that
crossings of diagrams should be appropriately labeled, since
we are in a multi-braided situation. At the diagramatic level, many of the
proofs become very simple. However, in this study the considerations
will be performed in the standard-algebraic way. For the reasons of
completeness, all the proofs are included in the paper. 

Through sections~3--5 we shall exclusively deal with
a given multi-braided  quantum  group $G$. In Section~3 we begin  with 
formulations of braided counterparts of concepts of the left, right 
and bi-covariance \cite{W}. As in the standard theory \cite{W},  the  notion 
of left covariance will be formulated by requiring  a  possibility 
of defining a left action $\ell_\Gamma\colon\Gamma\rightarrow
\cal{A}\otimes \Gamma$ of $G$
on $\Gamma$. Similarly,  right 
covariance will be characterized  by  a  possibility  of  a  right 
action $\rac\colon\Gamma\rightarrow\Gamma\otimes\cal{A}$.
The notion of bicovariance is a symbiosis of the previous two.

   Our attention will be  then  confined  to  the  left-covariant 
structures.  As  we  shall  see,  left  covariance  implies   left 
$\sigma$-covariance and, consequently, left $\sigma_n$-covariance,  for 
each $n\in\Bbb{Z}$. 
The corresponding  flip  over  operators $\ltw{\sigma_n}$ naturally  describe 
twisting properties of the left action $\ell_\Gamma$. Besides  the  study  of 
properties of maps $\ell_\Gamma$ and $\ltw{\sigma_n}$, and their
interrelations,
we shall also analyze the internal structure of  left-covariant calculi.
It turns out that the situation is more or less the same as
in the standard theory \cite{W}.
As a left/right $\cal{A}$-module, every left-covariant $\Gamma$ is free  and 
can  be  invariantly  decomposed   as
$$\Gamma\leftrightarrow\cal{A}\otimes\Gamma_{\inv}\leftrightarrow\Gamma_{\inv}
\otimes\cal{A}$$
where $\Gamma_{\inv}$ is 
the space of left-invariant elements of $\Gamma$. We shall also  prove  a 
braided generalization of the structure theorem \cite{W} by establishing 
a natural correspondence  between  (classes  of  isomorphic)
left-covariant $\Gamma$ and certain lineals $\cal{R}\subseteq\ker(\e)$.

    However, a full analogy with \cite{W} breaks, because $\cal{R}$
is generally not a right ideal in $\cal{A}$, but a right ideal in
a simplified \cite{D-BQG} algebra $\cal{A}_0$
obtained from $\cal{A}$ by an appropriate change of the  product. 
The lineal $\cal{R}$ should also be  left-invariant
with respect to the action of $\tau$.

    Concerning the concept of the right  covariance,  it  is  in  some 
sense symmetric to that  of the left  covariance.  For  this  reason, 
we shall not repeat completely analogous considerations for
right-covariant calculi.  The  most  important  properties  of   them  are 
collected,  without  proofs,  in  Appendix~A.   In   particular, 
right covariance implies right $\sigma_n$-covariance, for each $n\in\Bbb{Z}$.

    The study of bicovariant differential structures is  the  topic 
of Section~4.
In the bicovariant case the action maps $\ell_\Gamma$ and $\rac$, as well
as the flip-over maps $\ltw{\sigma_n},\rtw{\sigma_n}$ are mutually
compatible, in a natural manner.

 We  shall characterize  bicovariance in terms of the 
corresponding right $\cal{A}_0$-ideals $\cal{R}$.  It  turns  out  that  the 
calculus is bicovariant if and only if $\cal{R}$ satisfies two  additional 
conditions. The first one correspond to the adjoint invariance  in  the 
standard theory \cite{W}. In its formulation, a braided analogue of the 
adjoint action of $G$ on itself appears naturally. For this  reason, 
the most  important  properties  of  this  map  are  collected  in 
Appendix~B. The second additional condition for $\cal{R}$, trivial in  the 
standard theory, consists in its  right $\tau$-invariance.

    In  Section~5,  we  shall  analyze   differential   structures 
which are covariant with respect to the  antipodal  map 
$\k\colon\cal{A}\rightarrow\cal{A}$. 
Such structures will be called $\k$-covariant.
As we  shall  see,  a  symbiosis  of  $\k$  and  left  covariance  is 
equivalent to bicovariance.

   In Section~6 we shall first introduce the concept
of a *-structure on a multi-braided quantum group. Then, we pass to
the study of *-covariant calculi, in the context of multi-braided quantum
groups.

    Besides other results we shall obtain
a characterization  of  *-covariant
left covariant structures $\Gamma$, in terms  of  the  corresponding  right 
$\cal{A}_0$-ideal $\cal{R}$. It turns out that a
left-covariant calculus is *-covariant iff
$\k(\cal{R})^*\subseteq\cal{R}$, which is identical as in the
standard theory \cite{W}.

    In this paper only  the  abstract 
theory will be presented. Concrete examples will be  included 
in the next part of the study,  after  developing  a
higher-order differential calculus. This will include  differential
structures over already considered groups, as well as new 
examples of `differential' multi-braided quantum groups coming from  the 
developed theory. Finally, let us mention that we shall assume here
{\it trivial} braiding properties of the differential
$d\colon\cal{A}\rightarrow\Gamma$.
Our philosophy is that the non-trivial braidings involving the differential
map should be interpreted as an extra structure given over the whole
differential calculus. 

\subsection*{Acknowledgment} I would like to thank Professor Zbigniew
Oziewicz for carefully reading old ChiWriter version of the paper, and
for various valuable remarks and observations. 

\section{The Concept of Braided Covariance}

 Let $\cal{A}$ be a  complex  unital associative  algebra. 
Let us denote by 
$m\colon\cal{A}\otimes\cal{A}\rightarrow\cal{A}$ the multiplication in
$\cal{A}$. The
algebra  $\cal{A}$  will  be  interpreted  as  consisting  of  smooth 
functions over a quantum  space  $X$.

By definition, a first order differential calculus over $X$  is  a 
unital $\cal{A}$-bimodule $\Gamma$, equipped with a linear map
$d\colon\cal{A}\rightarrow\Gamma$ satisfying the Leibniz rule
\begin{equation}\label{21}
dm=\mlG(\id\otimes d)+\mrG(d\otimes\id),
\end{equation}
and such that $\ilG=\mlG(\id\otimes d)\colon\cal{A}\otimes 
\cal{A}\rightarrow\Gamma$ is surjective.
Here, $\mlG\colon\cal{A}\otimes\Gamma\rightarrow\Gamma$ and
$\mrG\colon\Gamma\otimes\cal{A}\rightarrow\Gamma$ are the left and the
right $\cal{A}$-module structures of $\Gamma$.

Let us observe that $d(1)=0$, and that the surjectivity of
$\ilG$ is equivalent to the surjectivity of $\irG\colon\cal{A}\otimes\cal{A}
\rightarrow\Gamma$, which is given by $\irG=\mrG(d\otimes\id)$. 

Now, let us assume that $X$ is a {\it braided} quantum space. In other words,
we have in addition a bijective braid operator
$\sigma\colon\cal{A}\otimes\cal{A}\rightarrow\cal{A}\otimes\cal{A}$
such that the following identities hold:
\begin{align}
(\id\otimes m)(\sigma\otimes \id)(\id\otimes\sigma)&=\sigma(m\otimes\id) 
\label{29}\\ 
(m\otimes\id)(\id\otimes\sigma)(\sigma\otimes\id)&=\sigma(\id\otimes  m). 
\label{210}
\end{align}
The operator $\sigma$ naturally induces a  structure  of  an  associative 
algebra  on $\cal{A}\otimes\cal{A}$ with the unit element $1\otimes 1\in 
\cal{A}\otimes \cal{A}$. Explicitly, the product is given by
$$ (a\otimes b)(q\otimes d)=a\sigma(b\otimes q)d. $$

We are going to analyze natural compatibility conditions between $\Gamma$
and $\sigma$.
\begin{defn}\label{def:24}
A first-order differential calculus $\Gamma$ over $X$ is called
left  $\sigma$-covariant  iff  there  exists  a  linear
operator $\ltw{\sigma}\colon\Gamma\otimes\cal{A}
\rightarrow\cal{A}\otimes\Gamma$ satisfying
\begin{equation}
\ltw{\sigma}(\ilG\otimes\id)=(\id\otimes\ilG)(\sigma\otimes 
\id)(\id\otimes\sigma).\label{216}
\end{equation}

Similarly, we say that $\Gamma$ is right $\sigma$-covariant  iff  there
exists a linear operator
$\rtw{\sigma}\colon\cal{A}\otimes
\Gamma\rightarrow\Gamma\otimes\cal{A}$ such that
\begin{equation}\label{217}
\rtw{\sigma}(\id\otimes\irG)=(\irG\otimes\id)(\id\otimes 
\sigma)(\sigma\otimes\id).
\end{equation}
Finally, $\Gamma$ is called bi-$\sigma$-covariant iff it  is  both 
right and left $\sigma$-covariant.
\end{defn}

The idea  beyond  this  definition  is  that  `twistings'
between elements from $\cal{A}$ and $\Gamma$ are performable `term  by 
term' such that twistings between the symbol $d$ and elements from 
$\cal{A}$ are trivial. 

It is easy to see that maps $\ltw{\sigma}$  and $\rtw{\sigma}$,
if they exist, are
uniquely determined by \eqref{216}  and \eqref{217} respectively.

Requirement \eqref{216} can be replaced by the equivalent
\begin{equation}
\ltw{\sigma}(\irG\otimes\id)=(\id\otimes\irG)
(\sigma\otimes\id)(\id\otimes \sigma).
\end{equation}

Similarly, the operator $\rtw{\sigma}$ appearing in the context of the
right $\sigma$-covariance can be characterized by
\begin{equation}
\rtw{\sigma}(\id\otimes\ilG)=(\ilG\otimes\id)(\id\otimes 
\sigma)(\sigma\otimes \id).
\end{equation}

In  the
following proposition, the most important  general  properties  of 
$\sigma$-covariant structures are collected.

\begin{pro}\label{pro:22}\bla{i} If $\Gamma$ is a left $\sigma$-covariant 
calculus, then
\begin{gather}
\ltw{\sigma}(\vartheta\otimes 1)=1\otimes\vartheta\qquad
\ltw{\sigma}(d\otimes
\id)=(\id\otimes d)\sigma\label{218}\\
(\id\otimes\ltw{\sigma})(\ltw{\sigma}\otimes\id)(
\id\otimes \sigma)=(\sigma\otimes \id)(\id\otimes
\ltw{\sigma})(\ltw{\sigma}\otimes \id).\label{220}
\end{gather}
The map $\ltw{\sigma}$ is surjective.   Its  kernel    is    an 
$\cal{A}$-subbimodule of $\Gamma\otimes\cal{A}$. Furthermore, we have
\begin{gather}
(\id\otimes\mlG)(\sigma\otimes\id)(\id\otimes\ltw{\sigma})=\ltw{\sigma}
(\mlG\otimes\id)\label{221}\\
(\id\otimes\mrG)(\ltw{\sigma}\otimes\id)(\id\otimes\sigma)=\ltw{\sigma}
(\mrG\otimes\id)\label{222}\\
(m\otimes\id)(\id\otimes\ltw{\sigma})(\ltw{\sigma}\otimes\id)=
\ltw{\sigma}(\id\otimes m)\label{223}
\end{gather}

\bla{ii} Similarly, if $\Gamma$ is right $\sigma$-covariant then
\begin{gather}
\rtw{\sigma}(1\otimes\vartheta)=\vartheta\otimes 1\qquad
\rtw{\sigma}(\id\otimes d)=(d\otimes \id)\sigma\label{226}\\
(\id\otimes\sigma)(\rtw{\sigma}\otimes\id)(\id\otimes 
\rtw{\sigma})=(\rtw{\sigma}\otimes \id)(\id\otimes\rtw{\sigma})
(\sigma\otimes \id).\label{227}
\end{gather}
The map $\rtw{\sigma}$ is surjective, its kernel is an 
$\cal{A}$-subbimodule of $\cal{A}\otimes \Gamma$ and the following
identities hold
\begin{gather}
(\id\otimes m)(\rtw{\sigma}\otimes\id)(\id\otimes\rtw{\sigma})=
\rtw{\sigma}(m\otimes\id)\label{228}\\
(\mrG\otimes\id)(\id\otimes\sigma)(\rtw{\sigma}\otimes\id)=
\rtw{\sigma}(\id\otimes \mrG)\label{229}\\
(\mlG\otimes\id)(\id\otimes\rtw{\sigma})(\sigma\otimes\id)=
\rtw{\sigma}(\id\otimes \mlG)\label{230}
\end{gather}

\bla{iii} Finally, if $\Gamma$ is bi-$\sigma$-covariant then
\begin{equation}\label{232}      
(\id\otimes\rtw{\sigma})(\sigma\otimes\id)(\id\otimes 
\ltw{\sigma})=(\ltw{\sigma}\otimes\id)(\id\otimes\sigma)(\rtw{\sigma}\otimes
\id).
\end{equation}
\end{pro}

\begin{pf}
Let us  assume  that $\Gamma$  is  left $\sigma$-covariant. 
Identities \eqref{218} are obvious. Let us check \eqref{221}--\eqref{223}.  
Using \eqref{216} we obtain
\begin{equation*}
\begin{split}
\ltw{\sigma}(\ilG\otimes m)&=(\id\otimes\ilG)(\sigma\otimes 
\id)(\id\otimes \sigma)(\id^2\otimes m)\\ 
&=(\id\otimes\ilG)(\sigma\otimes\id)(\id\otimes 
m\otimes \id)(\id^2\otimes \sigma)(\id\otimes \sigma\otimes \id)\\
&=(m\otimes\ilG)(\id\otimes\sigma\otimes 
\id)(\sigma\otimes\sigma)(\id\otimes\sigma\otimes \id)\\
&=(m\otimes\id)(\id\otimes\ltw{\sigma})(\id\otimes\ilG\otimes 
\id)(\sigma\otimes \id^2 )(\id\otimes \sigma\otimes \id)\\
&=(m\otimes \id)(\id\otimes \ltw{\sigma})(\ltw{\sigma}\otimes\id)
(\ilG\otimes\id^2).
\end{split}
\end{equation*}
Furthermore,
\begin{equation*}
\begin{split}
\ltw{\sigma}(\mlG\otimes\id)(\id\otimes\ilG\otimes\id)&=
\ltw{\sigma}(\ilG\otimes\id)(m\otimes\id^2)=(\id\otimes
\ilG)(\sigma\otimes \id)(m\otimes \sigma)\\ 
&=(\id\otimes \ilG)(\id\otimes m\otimes  \id)(\sigma\otimes 
\id^2)(\id\otimes \sigma\otimes \id)(\id^2\otimes \sigma)\\
&=(\id\otimes \mlG)(\sigma\otimes \id)(\id^2\otimes 
\ilG)(\id\otimes \sigma\otimes \id)(\id^2\otimes \sigma)\\
&=(\id\otimes\mlG)(\sigma\otimes\id)(\id\otimes 
\ltw{\sigma})(\id\otimes \ilG\otimes\id).
\end{split}
\end{equation*}
Similarly, we find
$$
\ltw{\sigma}(\mrG\otimes \id)(\ilG\otimes\id^2)=\ltw{\sigma}\Bigl\{\ilG
(\id\otimes m)\otimes\id\Bigr\}-\ltw{\sigma}\Bigl\{\ilG(m\otimes 
\id)\otimes \id\Bigr\}.
$$
The first term on the right-hand side of the above equality is equal to
$$
(\id\otimes\ilG)(\sigma\otimes\id)(\id\otimes 
\sigma)(\id\otimes m\otimes\id)=\bigl(\id\otimes\ilG 
(\id\otimes m)\bigr)(\sigma\otimes\id^2)(\id\otimes\sigma\otimes 
\id)(\id^2\otimes\sigma),
$$
while the second term reads
$$
(\id\otimes\ilG)(\sigma\otimes\id)(\id\otimes 
\sigma)(m\otimes \id^2)=(\id\otimes \ilG)(\id\otimes  m\otimes 
\id)(\sigma\otimes\id^2)(\id\otimes\sigma\otimes\id)(\id^2\otimes 
\sigma).
$$
Combining the last two expressions we conclude
\begin{multline*}
\ltw{\sigma}(\mrG\otimes\id)(\ilG\otimes\id^2)=(\id\otimes 
\mrG)(\id\otimes\ilG\otimes\id)(\sigma\otimes 
\id^2)(\id\otimes \sigma\otimes \id)(\id^2\otimes\sigma)\\
=(\id\otimes\mrG)(\ltw{\sigma}\otimes \id)(\ilG\otimes\sigma).
\end{multline*}

    We prove \eqref{220}. Direct transformations give
\begin{equation*}
\begin{split}
&(\id\otimes\ltw{\sigma})(\ltw{\sigma}\otimes\id)(\ilG\otimes\sigma)=\\
=&(\id\otimes\ltw{\sigma})(\id\otimes\ilG\otimes\id)(\sigma\otimes\id^2)
(\id\otimes\sigma\otimes \id)(\id^2\otimes\sigma)\\ 
=&(\id^2\otimes\ilG)(\id\otimes\sigma\otimes 
\id)(\sigma\otimes\sigma)(\id\otimes\sigma\otimes\id)(\id^2\otimes \sigma)\\ 
=&(\id^2\otimes \ilG)(\id\otimes \sigma\otimes 
\id)(\sigma\otimes\id^2)(\id\otimes\sigma\otimes\id)(\id^2\otimes 
\sigma)(\id\otimes\sigma\otimes\id)\\
=&(\id^2\otimes\ilG)(\sigma\otimes\id^2)(\id\otimes 
\sigma\otimes\id)(\sigma\otimes\sigma)(\id\otimes\sigma\otimes\id)\\
=&(\sigma\otimes\id)(\id\otimes\ltw{\sigma})(\id\otimes\ilG 
\otimes\id)(\sigma\otimes \id^2)(\id\otimes\sigma\otimes\id)\\
=&(\sigma\otimes\id)(\id\otimes\ltw{\sigma})(\ltw{\sigma}\otimes\id)
(\ilG\otimes\id).
\end{split}
\end{equation*}

 To prove the surjectivity of $\ltw{\sigma}$, it is  sufficient  to  check
that the elements of the form $a\otimes bd(q)$ belong   to
$\im(\ltw{\sigma})$. Let us define
$$
\omega=(\mlG\otimes\id)(\id\otimes d\otimes\id)(\id\otimes 
\sigma^{-1})(\sigma^{-1}\otimes\id)(a\otimes b\otimes q).
$$
Using \eqref{218} and \eqref{221} we obtain
$$
\ltw{\sigma}(\omega)
=(\id\otimes\mlG)(\sigma\otimes\id)(\id\otimes\ltw{\sigma})(\id\otimes 
d\otimes\id)(\id\otimes\sigma^{-1})
(\sigma^{-1}\otimes\id)(a\otimes b\otimes q)= 
a\otimes   bd(q).
$$
The fact that $\ker(\ltw{\sigma})$ is an
$\cal{A}$-subbimodule of $\Gamma\otimes \cal{A}$
directly follows from equalities \eqref{221} and \eqref{223}.

 In such a way we have shown \bla{i}. The right $\sigma$-covariance  case 
can be treated in a similar  manner.  Finally,  if $\Gamma$  is 
bi-$\sigma$-covariant then
\begin{equation*}
\begin{split}
&(\id\otimes\rtw{\sigma})(\sigma\otimes\id)(\id\otimes 
\ltw{\sigma})(\id\otimes\ilG\otimes\id)=\\
=&(\id\otimes\rtw{\sigma})(\sigma\otimes\id)(\id^2\otimes\ilG)
(\id\otimes \sigma\otimes \id)(\id^2\otimes \sigma)\\
=&(\id\otimes \ilG\otimes \id)(\id^2\otimes\sigma)(\id\otimes 
\sigma\otimes\id)(\sigma\otimes\id^2)(\id\otimes\sigma\otimes 
\id)(\id^2\otimes \sigma)\\ 
=&(\id\otimes\ilG\otimes\id)(\sigma\otimes 
\sigma)(\id\otimes \sigma\otimes \id)(\sigma\otimes \sigma)\\
=&(\id\otimes \ilG\otimes \id)(\sigma\otimes  \id^2)(\id\otimes 
\sigma\otimes\id)(\id^2\otimes   \sigma)(\id\otimes\sigma\otimes 
\id)(\sigma\otimes \id^2)\\ 
=&(\ltw{\sigma}\otimes\id)(\ilG\otimes\sigma)(\id\otimes 
\sigma\otimes \id)(\sigma\otimes \id^2)\\
=&(\ltw{\sigma}\otimes \id)(\id\otimes
\sigma)(\rtw{\sigma}\otimes \id)(\id\otimes \ilG\otimes \id),
\end{split}
\end{equation*}
and this completes the proof. 
\end{pf}

 It is easy to construct `pathological' examples in which maps 
$\ltw{\sigma}$ or $\rtw{\sigma}$ are not injective. However,
besides certain technical complications such a structure gives nothing
essentially new. For this reason, we shall assume from this moment that
every left/right $\sigma$-covariant  calculus  we  are  dealing   with 
possesses bijective flip-over operator $\ltw{\sigma}$ or $\rtw{\sigma}$.
Modulo this assumption, left $\sigma$-covariance and right 
$\sigma^{-1}$-covariance are equivalent properties. In other words, 
\begin{equation}
(\ltw{\sigma})^{-1}=\rtw{(\sigma^{-1})}\qquad
(\rtw{\sigma})^{-1}=\ltw{(\sigma^{-1})}.
\end{equation}

    Now,  we  shall  generalize  the  previous  consideration   to 
situations  in  which,  instead  of  one,  a  system  of   mutually 
compatible braided quantum space \cite{D-BQG} structures on $X$
appears.

\begin{defn}\label{def:25}
Let us assume that $\cal{A}$ is equipped with a braid system $\cal{T}$.
Then we shall say that $X$ is a $\cal{T}$-braided quantum space.
\end{defn}

\begin{defn}\label{def:26} A first order calculus $\Gamma$
over a $\cal{T}$-braided quantum space $X$ is called
left/right/bi $\cal{T}$-covariant iff it is left/right/bi
$\gamma$-covariant for each braiding $\gamma\in\cal{T}$.
\end{defn}

As explained in \cite{D-BQG}-Appendix, every braid system $\cal{T}$
can be naturally {\it completed}. The completed system $\cal{T}^*$ is
defined as the minimal extension of $\cal{T}$, invariant under
ternar operations of the form $\delta=\alpha\beta^{-1}\gamma$. Explicitly,
$\cal{T}^*$ is the union of systems $\cal{T}_n$, where $\cal{T}_0=\cal{T}$
and $\cal{T}_{n+1}$ is obtained from $\cal{T}_n$ by applying the above
mentioned operations. 

\begin{pro}\label{pro:23} Let $X$ be a $\cal{T}$-braided quantum space
and $\Gamma$ a first-order calculus over $X$. 

\smallskip
\bla{i} If $\Gamma$ is left $\cal{T}$-covariant then
it is also left $\cal{T}^*$-covariant. We have
\begin{align}
\ltw{(\alpha\beta^{-1}\gamma)}&=\ltw{\alpha}(\ltw{\beta})^{-1}\ltw{\gamma}
\label{234}\\
(\id\otimes\ltw{\alpha})(\ltw{\beta}\otimes\id)(\id\otimes\gamma 
)&=(\gamma\otimes\id)(\id\otimes\ltw{\beta})(\ltw{\alpha}\otimes\id), 
\label{235}
\end{align}
for each $\alpha,\beta,\gamma \in\cal{T}^*$.

\smallskip
\bla{ii} Similarly, if $\Gamma$ is right $\cal{T}$-covariant then it 
is right $\cal{T}^*$-covariant, too. We have
\begin{align}
\rtw{(\alpha\beta^{-1}\gamma)}&= \rtw{\alpha}(\rtw{\beta})^{-1}
\rtw{\gamma}\label{236}\\
(\id\otimes\alpha)(\rtw{\beta}\otimes\id)(\id\otimes\rtw{\gamma})
&=(\rtw{\gamma}\otimes\id)(\id\otimes\rtw{\beta})(\alpha\otimes\id) 
\label{237}
\end{align}
for each $\alpha,\beta,\gamma\in\cal{T}^*$.

\smallskip
\bla{iii} If $\Gamma$ is bi-$\cal{T}$-covariant, it  is  consequently 
also bi-$\cal{T}^*$-covariant and
\begin{equation}\label{238}
(\id\otimes\rtw{\alpha})(\beta\otimes\id)(\id\otimes\ltw{\gamma}
)=(\ltw{\gamma}\otimes\id)(\id\otimes\beta)(\rtw{\alpha}\otimes\id) 
\end{equation}
for each $\alpha,\beta,\gamma\in\cal{T}^*$.
\end{pro}

\begin{pf}
Let us assume that $\Gamma$ is left $\cal{T}$-covariant. Then,
\begin{multline*}
\ltw{\alpha}(\ltw{\beta})^{-1}\ltw{\gamma}(\ilG\otimes
\id)=(\id\otimes\ilG)(\alpha\otimes\id)(\id\otimes
\alpha\beta^{-1})(\beta^{-1}\gamma\otimes \id)(\id\otimes \gamma)\\ 
=(\id\otimes\ilG)(\alpha\beta^{-1}\gamma\otimes 
\id)(\id\otimes \alpha \beta^{-1}\gamma),
\end{multline*}
for  each  $\alpha,\beta,\gamma\in\cal{T}$. 
This means that $\Gamma$  is  left
$\alpha\beta^{-1}\gamma$-covariant and \eqref{234}
holds for the braidings from $\cal{T}$.
Now, we can proceed inductively and conclude that $\Gamma$ is 
left $\cal{T}_n$-covariant for each $n\in\Bbb{N}$, and that
\eqref{234} holds on $\cal{T}^*$.

Similarly, if $\Gamma$ is right $\cal{T}$-covariant then
\begin{multline*}
\rtw{\alpha}(\rtw{\beta})^{-1}\rtw{\gamma}(\id\otimes\irG)=(\irG
\otimes\id)(\id\otimes\alpha)(\alpha\beta^{-1}\otimes 
\id)(\id\otimes \beta^{-1}\gamma)(\gamma \otimes \id)\\
=(\irG \otimes \id)(\id\otimes \alpha \beta^{-1}\gamma)(\alpha\beta^{-1} 
\gamma\otimes\id),
\end{multline*}
for each $\alpha,\beta,\gamma\in\cal{T}$.
This implies that $\Gamma$ is right $\cal{T}^*$-covariant and  that 
\eqref{236} holds for each $\alpha,\beta,\gamma\in\cal{T}^*$.

    Identities \eqref{235} and \eqref{237}--\eqref{238}
can be derived in
essentially the  same  manner  as  it  is  done  in  the  proof  of 
Proposition~\ref{pro:22}, in the case of a single flip-over operator.
\end{pf}

    From this moment, as well as through the next  three  sections we
shall deal exclusively with braided quantum  groups, in the sense of
\cite{D-BQG}. Let $G$ be such a group, represented by $\cal{A}$. 
We shall denote by $\phi\colon\cal{A}\rightarrow\cal{A}
\otimes\cal{A}$ the coproduct map, and by
$\e\colon\cal{A}\rightarrow\Bbb{C}$ and
$\k\colon\cal{A}\rightarrow\cal{A}$ the counit and the antipode map
respectively. Let $\sigma\colon\cal{A}\otimes\cal{A}\rightarrow\cal{A}
\otimes\cal{A}$ be the intrinsic braid operator. 

  As explained in \cite{D-BQG}, twisting properties of the coproduct and
the antipode are not properly
expressible in terms a single braid operator $\sigma$.
This is the place where a `secondary' braid operator naturally
enters the game. Explicitly, it is given by
\begin{equation}
\tau=(\e\otimes\id^2)(\sigma^{-1}\otimes\id)(\id\otimes\phi)\sigma=
(\id^2\otimes\e)(\id\otimes\sigma^{-1})(\phi\otimes\id)\sigma. 
\end{equation}

The operators $\{\sigma,\tau\}$ form a braid system, and the completion 
$\cal{F}=\{\sigma,\tau\}^*$ is consisting of maps of the form
$$
\sigma_n=(\sigma\tau^{-1})^{n-1}\sigma=\sigma(\tau^{-1}\sigma)^{n-1},
$$
where $n\in\Bbb{Z}$.

\begin{pro}\label{pro:24}
\bla{i} If $\Gamma$ is left $\sigma$-covariant then it is also left
$\cal{F}$-covariant and
\begin{gather}
\ltw{\tau}=(\id^2\otimes\e)(\id\otimes(\ltw{\sigma})^{-1})
(\phi \otimes \id)\ltw{\sigma}\label{239}\\
(\ltw{\tau})^{-1}=(\e\otimes\id^2)(\ltw{\sigma}\otimes\id)(\id\otimes\phi)
(\ltw{\sigma})^{-1}\label{240}\\
(\e\otimes \id)\ltw{\tau}= \id\otimes\e.\label{241}                            
\end{gather}
Moreover, the following twisting properties hold
\begin{equation}
(\id\otimes\ltw{\sigma_n})(\ltw{\sigma_m}\otimes\id)(\id\otimes\phi)
=(\phi\otimes\id)\ltw{\sigma_{m+n}}\label{242}
\end{equation}
for each $n,m\in\Bbb{Z}$.

\bla{ii} Similarly, if $\Gamma$ is right $\sigma$-covariant then it is 
also right $\cal{F}$-covariant and
\begin{gather}
\rtw{\tau}=(\e\otimes \id^2)((\rtw{\sigma})^{-1}\otimes \id)(\id\otimes
\phi)\rtw{\sigma}\label{243}\\
(\rtw{\tau})^{-1}=(\id^2\otimes\e)(\id\otimes\rtw{\sigma})(\phi \otimes\id)
(\rtw{\sigma})^{-1}\label{244}\\
(\id\otimes\e)\rtw{\tau}=\e\otimes\id.\label{245}
\end{gather}
We also have
\begin{equation}
(\rtw{\sigma_n}\otimes\id)(\id\otimes\rtw{\sigma_m})(\phi\otimes\id)
=(\id\otimes\phi)\rtw{\sigma_{n+m}},           \label{246}
\end{equation}
for each $n,m\in\Bbb{Z}$.
\end{pro}

\begin{pf}
Let us assume left $\sigma$-covariance of $\Gamma$, and consider a
map $\xi\colon\Gamma\otimes \cal{A}\rightarrow\cal{A}\otimes \Gamma$
determined by the right
hand side of \eqref{239}. Direct transformations give
\begin{equation*}
\begin{split}
\xi(\ilG\otimes \id)&=(\id^2\otimes\e)(\id\otimes  (\ltw{\sigma})^{-1})
(\phi\otimes \id)\ltw{\sigma}(\ilG  \otimes \id)\\ 
&=(\id\otimes\ilG\otimes\e)(\id^2\otimes 
\sigma^{-1})(\id\otimes\sigma^{-1}\otimes\id)(\phi \otimes 
\id^2)(\sigma\otimes \id)(\id\otimes \sigma)\\
&=(\id\otimes\ilG\otimes\e)(\tau\otimes \sigma^{-1})(\id\otimes 
\phi \otimes \id)(\id\otimes \sigma)\\
&=(\id\otimes\ilG)(\tau\otimes\id)(\id\otimes\tau).
\end{split}
\end{equation*}
Consequently,  $\Gamma$ is  left  $\tau$-covariant  and $\xi=\ltw{\tau}$. 
According to Proposition~\ref{pro:23} the  calculus  is  left
$\cal{F}$-covariant.

   Let us denote by $\psi$ a map determined by the right hand  side  of 
\eqref{240}. We have then
\begin{equation*}
\begin{split}
\psi(\id\otimes\irG)&=(\e\otimes\id^2)(\ltw{\sigma}\otimes 
\id)(\id\otimes \phi)(\ltw{\sigma})^{-1}(\id\otimes\irG)\\
&=(\e\otimes\irG\otimes\id)(\sigma\otimes\id^2)(\id\otimes 
\sigma\otimes\id)(\id^2\otimes\phi)(\id\otimes 
\sigma^{-1})(\sigma^{-1}\otimes\id)\\
&=(\e\otimes\irG\otimes \id)(\sigma\otimes \tau^{-1})(\id\otimes  \phi 
\otimes\id)(\sigma^{-1}\otimes\id)\\
&=(\e\otimes\irG\otimes
\id)(\id^2\otimes \tau^{-1})(\id\otimes  \tau^{-1}\otimes\id)
(\phi\otimes\id)\\
&=( \irG\otimes\id)(\id\otimes\tau^{-1})(\tau^{-1}\otimes\id).
\end{split}
\end{equation*}
Consequently, $\Gamma$ is right $\tau^{-1}$-covariant 
and $\psi=\rtw{(\tau^{-1})}=(\ltw{\tau})^{-1}$.

   Let us prove the twisting property \eqref{242}. Using the standard braid
relations we obtain
\begin{multline*}
    (\id\otimes\ltw{\sigma_n})(\ltw{\sigma_m}\otimes\id)(\id\otimes 
\phi)(\ilG\otimes\id)=\\ 
=(\id\otimes\ltw{\sigma_n})(\id\otimes 
\ilG\otimes\id)(\sigma_m\otimes\id^2)(\id\otimes\sigma_m\otimes 
\id)(\id^2\otimes \phi)\\
=(\id^2\otimes \ilG)(\id\otimes\sigma_n\otimes \id)
(\sigma_m\otimes
\sigma_n)(\id\otimes\sigma_m\otimes \id)(\id^2\otimes\phi)\\
=(\id^2\otimes\ilG)(\id\otimes 
\sigma_n\otimes\id)(\sigma_m\otimes  \id^2)(\id\otimes  \phi\otimes 
\id)(\id\otimes \sigma_{m+n})\\
=(\phi\otimes\ilG)(\sigma_{n+m}\otimes\id)(\id\otimes 
\sigma_{n+m}) 
=(\phi \otimes \id)\ltw{\sigma_{n+m}}(\ilG\otimes \id).
\end{multline*}
The case \bla{ii}, when $\Gamma$ is right $\sigma$-covariant, can  be
treated in a similar way.
\end{pf}

    Finally, let us describe  twisting relations  between   the 
antipode $\k$ and a $\sigma$-covariant calculus $\Gamma$.

\begin{pro}\label{pro:25}
If $\Gamma$ is left $\sigma$-covariant then
\begin{equation}\label{247}
\ltw{\sigma_n}(\id\otimes\k)=(\k\otimes \id)\ltw{\sigma_{-n}}.
\end{equation}
Similarly, if $\Gamma$ is right $\cal{F}$-covariant then
\begin{equation}\label{248}
\rtw{\sigma_n}(\k\otimes \id)=(\id\otimes\k)\rtw{\sigma_{-n}},
\end{equation}
for  each $n\in\Bbb{Z}$.
\end{pro}

\begin{pf}
Let us assume that $\Gamma$  is  left $\sigma$-covariant. We have
\begin{multline*}
\ltw{\sigma}(\ilG\otimes\k)=(\id\otimes\ilG)(\sigma_n\otimes 
\id)(\id\otimes\sigma_n)(\id\otimes\k)=(\k\otimes\ilG
)(\sigma_{-n}\otimes \id)(\id\otimes\sigma_{-n})\\ 
=(\k\otimes \id)\ltw{\sigma_{-n}}(\ilG\otimes \id).
\end{multline*}
If the calculus is right $\cal{F}$-covariant then 
\begin{multline*}   
\rtw{\sigma_n}(\k\otimes 
\irG)=(\irG\otimes\id)(\id\otimes\sigma_n)(\sigma_n\otimes 
\id)(\k\otimes\id^2)=(\irG\otimes\k)(\id\otimes 
\sigma_{-n})(\sigma_{-n}\otimes \id)\\ 
=(\id\otimes \k) \rtw{\sigma_{-n}}(\id\otimes\irG),
\end{multline*}
for each $n\in\Bbb{Z}$.
\end{pf}

\section{The Structure of Left-Covariant Calculi}

We pass to definitions of  first  order  differential  structures 
which are covariant with respect to the comultiplication map $\phi 
\colon\cal{A}\rightarrow\cal{A}\otimes\cal{A}$.

\begin{defn}\label{def:31}
A first-order calculus $\Gamma$  over $G$
is  called  right-covariant
iff there exists  a  linear  map $\rac\colon\Gamma\rightarrow\Gamma 
\otimes \cal{A}$ such that
\begin{equation}\label{31}
\rac\irG=( \irG\otimes  m)(\id\otimes\sigma  \otimes\id)(\phi 
\otimes\phi).
\end{equation}
The map $\rac$ is called the right action of $G$ on $\Gamma$. It is
uniquely determined by the above condition.
\end{defn}

\begin{defn}\label{def:32}
The calculus $\Gamma$ is called left-covariant 
iff there exists a left action map $\lac\colon\Gamma\rightarrow\cal{A} 
\otimes \Gamma$ satisfying
\begin{equation}\label{32}
\lac\ilG=(m\otimes\ilG)(\id\otimes\sigma\otimes 
\id)(\phi \otimes \phi).
\end{equation}
The map $\lac$ is uniquely determined by this condition. 
\end{defn}

\begin{defn}\label{def:33}
We shall say that the calculus $\Gamma$ is bicovariant,  iff  it 
is both left and right-covariant.
\end{defn}

The above definitions naturally formulate braided generalizations
of standard concepts of right/left and bi-covariance in the standard theory 
\cite{W}. Throughout the rest of the section, we shall consider
left-covariant differential structures.

\begin{pro}\label{pro:31} We have
\begin{align}
\lac d&=(\id\otimes d)\phi\label{33}\\
\lac\mlG&=(m\otimes\mlG)(\id\otimes\sigma\otimes\id)(\phi 
\otimes \lac).\label{34}
\end{align}
\end{pro}

\begin{pf} Identity \eqref{33} is a direct consequence of \eqref{32}.  
To prove \eqref{34}, we start from \eqref{32} and apply elementary properties
of the product and the coproduct maps:
\begin{equation*}
\begin{split}
(\lac\mlG&)(\id\otimes\ilG)=\lac\ilG(m\otimes  \id)=(m\otimes 
\ilG)(\id\otimes \sigma \otimes \id)(\phi \otimes \phi)(m\otimes\id)\\
&=(m\otimes\ilG )(m\otimes\sigma\otimes\id)(\id^2\otimes 
m\otimes \id^2)(\id\otimes\sigma\otimes\id^3)(\phi  \otimes\phi 
\otimes \phi)\\ 
&=(m\otimes  \ilG)(m\otimes  \id\otimes  m\otimes  \id)(\id^2\otimes 
\sigma  \otimes  \id^2)(\id\otimes  \sigma  \otimes  \sigma   \otimes 
\id)(\phi \otimes \phi \otimes \phi )\\ 
&=(m\otimes  \ilG)(\id^2\otimes  m\otimes   \id)\bigl[\id\otimes   \sigma 
(\id\otimes  m)\otimes  \id^2\bigr](\id^3\otimes  \sigma  \otimes   \id)(\phi 
\otimes \phi \otimes \phi )\\ 
&=(m\otimes\mlG)\bigl[\id\otimes \sigma 
(\id\otimes  m)\otimes  \ilG\bigr](\id^3\otimes\sigma\otimes\id)(\phi 
\otimes \phi \otimes \phi )\\ 
&=(m\otimes \mlG )(\id\otimes \sigma \otimes \id)(\phi 
\otimes \lac\ilG). \qed
\end{split}
\end{equation*}
\renewcommand{\qed}{}
\end{pf}

It is worth noticing that
\begin{equation}\label{35}
\lac\irG=(m\otimes\irG)(\id\otimes\sigma\otimes\id)(\phi 
\otimes\phi),
\end{equation}
which also characterizes the map $\lac$. The following proposition
shows that $\lac$ gives a left $\cal{A}$-comodule structure on $\Gamma$. 

\begin{pro}\label{pro:32} We have
\begin{gather}
(\e\otimes \id)\lac=\id\label{36}\\
(\phi \otimes \id)\lac=(\id\otimes \lac)\lac.\label{37}
\end{gather}
\end{pro}

\begin{pf} Applying \eqref{32} and performing further elementary
transormations with the counit we obtain
\begin{equation*}
\begin{split}
(\e\otimes \id)\lac\ilG&=(\e\otimes\id)(m\otimes \ilG)(\id\otimes 
\sigma \otimes\id)(\phi\otimes\phi)\\
&=(\e\otimes\e\otimes\ilG)(\sigma^{-1}\tau \otimes\id^2)(\id\otimes 
\sigma\otimes\id)(\phi\otimes\phi)\\ 
&=(\e\otimes\e\otimes \ilG)(\sigma^{-1}\otimes \id^2)(\id\otimes \phi 
\otimes\id)(\sigma\otimes\id)(\id\otimes\phi 
)\\
&=(\e\otimes\e\otimes \ilG)(\id\otimes \tau \otimes \id)(\phi
\otimes\phi)=\ilG.
\end{split}
\end{equation*}
Furthermore,
\begin{equation*}
\begin{split}
(&\id\otimes\lac)\lac\ilG=(m\otimes m\otimes 
\ilG)(\id^3\otimes\sigma\otimes\id)(\id^2\otimes  \phi 
\otimes \phi)(\id\otimes \sigma \otimes \id)(\phi \otimes 
\phi)\\
&=(m\otimes m\otimes   \ilG)(\id^3\otimes\sigma\otimes 
\id)(\id^2\otimes  \phi  \otimes  \id^2)(\id\otimes\sigma\otimes 
\id^2)\bigl(\phi \otimes (\id\otimes \phi )\phi\bigr)\\ 
&=(m\otimes m\otimes    \ilG)(\id^3 \otimes    \sigma    \otimes 
\id)(\id^2 \otimes  \phi   \otimes   \id^2)(\id\otimes\sigma\otimes 
\id^2)(\phi \otimes \phi \otimes \id)(\id\otimes \phi )\\
&=(m\otimes m\otimes\ilG)(\id\otimes    \sigma\otimes 
\id^3)(\id^2\otimes  \phi  \otimes\id^2)(\id^2\otimes\sigma\otimes 
\id)(\id\otimes \phi \otimes \id^2)(\phi \otimes \phi )\\
&=(m\otimes m\otimes \ilG)(\id\otimes  \sigma  \otimes  \id^3)(\phi 
\otimes \phi \otimes \id^2)(\id\otimes \sigma \otimes \id)(\phi \otimes 
\phi )\\
&=(\phi m\otimes \ilG )(\id\otimes \sigma  \otimes  \id)(\phi  \otimes 
\phi )=(\phi \otimes \id)\lac\ilG,
\end{split}
\end{equation*}
which completes the proof. We have used the `octagonal' compatibility
property between $\phi$ and $\sigma$. 
\end{pf}

As we shall now see, every  left-covariant  differential  calculus
$\Gamma$ is left $\sigma$-covariant.  According  to
Proposition~\ref{pro:24}, this means that $\Gamma$ is left
$\cal{F}$-covariant, too.

\begin{pro}\label{pro:33} \bla{i}  The  calculus $\Gamma$ is,   being 
left-covariant, also left $\cal{F}$-covariant. 

\smallskip
\bla{ii} The diagram
\begin{equation}\label{38}
\begin{CD}
\cal{A}\otimes\Gamma\otimes\cal{A}\otimes\cal{A} @>{\mbox{$\k\otimes\lac\mrG
\otimes\k$}}>> \cal{A}\otimes\cal{A}\otimes\Gamma\otimes\cal{A}\\
@A{\mbox{$\lac\otimes\phi$}}AA  @VV{\mbox{$m\otimes\mrG$}}V\\
\Gamma\otimes\cal{A}
@>>{\mbox{$\quad\ltw{\sigma}\quad$}}>
\cal{A}\otimes\Gamma
\end{CD}
\end{equation}
is commutative.
\end{pro}

\begin{pf}
Let $\xi\colon\Gamma \otimes \cal{A}\rightarrow\cal{A}\otimes\Gamma$
be a map determined by
$$
 \xi=(m\otimes\mrG)(\k\otimes \lac\mrG\otimes \k)(\lac\otimes \phi).
$$
We shall prove that $\xi$ satisfies a characteristic  property  for 
the flip-over operator  $\ltw{\sigma}$. A direct computation gives
\begin{multline*}
\xi(\irG\otimes  \id)=(m\otimes\mrG)(\k\otimes\lac\mrG\otimes 
\k)(m\otimes  \irG\otimes \id^2)(\id\otimes \sigma \otimes \id^3)(\phi 
\otimes \phi \otimes \phi )\\ 
=(m\otimes\mrG)(m\otimes \lac\irG\otimes\k)(\k\otimes  \k\otimes 
\id\otimes  m\otimes  \id)(\sigma_{-2}\otimes\id^4)(\id\otimes  \sigma 
\otimes \id^3)(\phi \otimes \phi \otimes \phi )\\
=(m\otimes\mrG)(m\otimes \lac\irG\otimes \id)(\k\otimes \k\otimes \id^3)A\\ 
=(m\otimes\mrG)(m\otimes  m\otimes   \irG\otimes   \id)(\id^3\otimes 
\sigma \otimes \id^2)(\k\otimes \k\otimes \phi  \otimes  \phi  \otimes 
\id)A\\
=(m\otimes   \irG)(m\otimes   \id^2\otimes   m)(\id^2\otimes   \sigma 
\otimes \id^2)\bigl[\k\otimes m(\k\otimes \id)\phi  \otimes  \id\otimes  \phi 
\otimes\id\bigr]A\\
=(m\otimes\irG)(m\otimes \sigma\otimes m)(\k\otimes 
  1\e\otimes \id\otimes \phi  \otimes  \id)A,
\end{multline*}
where we have introduced $A=(\id\otimes \phi \otimes \id^2)(\sigma_{-1}\otimes 
m\otimes\k)(\id\otimes \phi \otimes \phi )$.

   The last term in the above sequence of
transformations can be further written as follows:
\begin{multline*}
(m\otimes\irG)(\k\otimes    \sigma    \otimes 
m)(\id^2\otimes  \phi  \otimes  \id)(\sigma_{-1}\otimes 
m\otimes\k)(\id\otimes \phi \otimes \phi)=\\
=(m\otimes    \irG)(\id\otimes   \sigma    \otimes 
m)(\sigma  \otimes  \phi   \otimes   \id)(\id\otimes 
\k\otimes m\otimes \k)(\id\otimes \phi  \otimes  \phi)\\
=(\id\otimes\irG)(\sigma \otimes \id)\bigl[\id\otimes
(m\otimes m)(\k\otimes \phi
m\otimes\k)(\phi\otimes\phi)\bigr]=(\id\otimes\irG)(\sigma \otimes 
\id)(\id\otimes \sigma).
\end{multline*}

   Thus, $\Gamma$ is left $\sigma$-covariant, $\chi=\ltw{\sigma}$ and 
diagram~\eqref{38} is 
commutative.  According  to  \bla{i}-Proposition~\ref{pro:24} the
calculus $\Gamma$ is automatically left $\cal{F}$-covariant.
\end{pf}

 The operator $\ltw{\sigma}$ figures in the right multiplicativity law for
the left action map.
\begin{pro}\label{pro:34} The diagram
\begin{equation}\label{39}
\begin{CD}
\cal{A}\otimes\Gamma\otimes\cal{A}\otimes\cal{A}
@>{\mbox{$\id\otimes\ltw{\sigma}
\otimes\id$}}>> \cal{A}\otimes\cal{A}\otimes\Gamma\otimes\cal{A}\\
@A{\mbox{$\lac\otimes\phi$}}AA  @VV{\mbox{$m\otimes\mrG$}}V\\
\Gamma\otimes\cal{A}
@>>{\mbox{$\quad\lac\mrG\quad$}}> \cal{A}\otimes\Gamma
\end{CD}
\end{equation}
is commutative.
\end{pro}

\begin{pf} According to Proposition~\ref{pro:32} and diagram~\eqref{38},
\begin{multline*}
(m\otimes\mrG)(\id\otimes \ltw{\sigma}\otimes\id)(\lac\otimes  \phi)=\\
=(m\otimes\mrG)\bigl[\id\otimes (m\otimes\mrG)(\k\otimes  \lac\mrG\otimes
\k)(\lac\otimes \phi )\otimes \id\bigr](\lac\otimes \phi )\\=(m\otimes 
\mrG)(m\otimes  \id^2\otimes  m)(\id\otimes\k\otimes   \lac
\mrG\otimes\k\otimes\id)(\phi  \otimes  \id^2\otimes  \phi 
)(\lac\otimes \phi )\\
=(m\otimes\mrG)(1\e\otimes \lac\mrG\otimes 1\e)(\lac
\otimes \phi)=\lac\mrG.\qed
\end{multline*}
\renewcommand{\qed}{}
\end{pf}

 Our next proposition describes twisting  properties  of  the  left
action map, with respect to the braid system $\cal{F}$.

\begin{pro}\label{pro:35} We have
\begin{equation}\label{310}
(\sigma_n\otimes\id)(\id\otimes\ltw{\sigma_m})(\lac\otimes\id)=
(\id\otimes\lac)\ltw{\sigma_{n+m}},
\end{equation}
for each $n,m\in\Bbb{Z}$. In particular, it follows that
\begin{equation}\label{311}
\ltw{\tau}=(\e\otimes  \id^2)(\sigma^{-1}\otimes   \id)(\id\otimes\lac 
)\ltw{\sigma}. 
\end{equation}
\end{pro}

\begin{pf} Using~\eqref{32} and the main properties of $\cal{F}$ we obtain 
\begin{multline*}
(\sigma_n\otimes   \id)(\id\otimes   \ltw{\sigma_m})(\lac\otimes
\id)(\ilG \otimes \id)= \\
=(\sigma_n  \otimes \id)(\id\otimes \ltw{\sigma_m})(m\otimes \ilG\otimes 
\id)(\id\otimes \sigma \otimes \id^2)(\phi \otimes \phi \otimes \id)\\ 
=(\sigma_n  \otimes \id)(m\otimes \id\otimes \ilG)(\id^2\otimes  \sigma_m 
\otimes \id)(\id\otimes \sigma \otimes \sigma_m)(\phi \otimes  \phi 
\otimes \id)\\ 
=(\id\otimes m\otimes \ilG)(\sigma_n  \otimes\id^3)(\id\otimes 
\sigma_n   \otimes  \id^2)(\id^2\otimes  \sigma_m   \otimes  \id)(\id\otimes 
\sigma \otimes \sigma_m)(\phi \otimes \phi \otimes \id)\\ 
=(\id\otimes m\otimes \ilG)(\sigma_n\otimes 
\sigma \otimes \id)(\id\otimes\sigma_m\otimes\id^2)(\id^2 \otimes 
\sigma_n\otimes \id)(\id^3 \otimes \sigma_m)(\phi \otimes \phi \otimes\id)\\ 
=(\id\otimes m\otimes  \ilG)(\sigma_n\otimes   \sigma   \otimes 
\id)(\id\otimes \sigma_m  \otimes \id^2)(\id^3\otimes \phi )(\phi  \otimes 
\sigma_{n+m})\\ 
=(\id\otimes m\otimes    \ilG)(\id^2\otimes    \sigma     \otimes 
\id)(\id\otimes \phi \otimes \phi )(\sigma_{n+m}\otimes  \id)(\id\otimes 
\sigma_{n+m})\\
=(\id\otimes  \lac\ilG)(\sigma_{n+m}\otimes  \id)(\id\otimes  \sigma_{n+m}
)=(\id\otimes \lac)\ltw{\sigma_{n+m}}(\ilG \otimes \id). \qed
\end{multline*}
\renewcommand{\qed}{}
\end{pf}

 We pass to the study of the internal  structure  of  left-covariant
calculi. For a given $\Gamma$, let $\Gamma_{\inv}$ be the  space  of 
left-invariant elements of $\Gamma$. In other words
\begin{equation}
\Gamma_{\inv}=\Bigl\{\vartheta\in\Gamma\mid\lac(\vartheta)=1
\otimes \vartheta\Bigr\}.\label{312}
\end{equation}

    Let $P\colon\Gamma\rightarrow\Gamma$  be a linear map defined by
\begin{equation}\label{313}
P=\mlG(\k\otimes \id)\lac.
\end{equation}
We are going to show that $P$ projects $\Gamma$ onto $\Gamma_{\inv}$.
Evidently, the elements of $\Gamma_{\inv}$ are $P$-invariant.
\begin{lem}\label{lem:36} We have
\begin{equation}\label{314}
 P\ilG=(\e\otimes Pd)\sigma^{-1}\tau.
\end{equation}
\end{lem}

\begin{pf}
Applying \eqref{32}--\eqref{33}, \eqref{313} and performing standard
braided transformations we obtain
\begin{equation*}
\begin{split}
P\ilG &=\mlG(\k\otimes  \id)(m\otimes   \ilG)(\id\otimes   \sigma 
\otimes \id)(\phi \otimes \phi )\\ 
&=\mlG(m\otimes \ilG)(\k\otimes \k\otimes  \id^2)(\tau\sigma^{-1}\tau 
\sigma^{-1}\tau\otimes  \id^2 )(\id\otimes  \sigma  \otimes  \id)(\phi 
\otimes \phi )\\ 
&=\mlG(m\otimes \ilG)(k\otimes k\otimes  \id^2)(\id\otimes 
\phi   \otimes   \id)(\tau\sigma^{-1}\tau\otimes\id)(\id\otimes \phi)\\
&=\mlG(\id\otimes   \ilG)(\k\otimes 1\e\otimes    \id)(\tau    \otimes 
\id)(\id\otimes   \phi   )\sigma^{-1}\tau\\
&=\mlG(\e\otimes\k\otimes d)(\id\otimes \phi)\sigma^{-1}\tau= 
(\e\otimes Pd)\sigma^{-1}\tau .\qed
\end{split}
\end{equation*}
\renewcommand{\qed}{}
\end{pf}
Now, it follows that  $P(\Gamma)\subseteq\Gamma_{\inv}$.  Indeed, according
to  the previous lemma,  it  is  sufficient  to  check  that
$Pd(\cal{A})\subseteq\Gamma_{\inv}$. We compute
\begin{equation*}
\begin{split}
\lac Pd&=\lac\mlG(\k\otimes  d)\phi=(m\otimes\mlG)(\id\otimes 
\sigma \otimes \id)(\phi \otimes \lac)(\k\otimes d)\phi\\ 
&=(m\otimes\mlG)(\id\otimes   \sigma   \otimes\id)(\sigma   \otimes 
\id^2)(\k\otimes \k\otimes \id\otimes d)(\phi \otimes \phi )\phi\\
&=(\id\otimes \mlG)(\sigma \otimes \id)(\id\otimes m\otimes \id)(\k\otimes 
\k\otimes \id\otimes  d)(\id\otimes  \phi  \otimes  \id)(\phi  \otimes 
\id)\phi\\ 
&=(\id\otimes \mlG)(\sigma \otimes \id)(\k\otimes 1\e\otimes d)
(\phi\otimes\id)\phi\\
&=1\otimes\mlG(\k\otimes d)\phi=1\otimes Pd.
\end{split}
\end{equation*}
Consequently, $P$ projects $\Gamma$  onto  $\Gamma_{\inv}$
and the composition
$$\pi=Pd=\mlG(\k\otimes   d)\phi\colon\cal{A}\rightarrow\Gamma_{\inv}$$     
is surjective.

 It is easy to see, by the  use  of \eqref{310},  that  the  flip-over 
operators $\ltw{\sigma_n}$ map $\Gamma_{\inv}\otimes \cal{A}$
onto $\cal{A}\otimes \Gamma_{\inv}$. Moreover, the corresponding restrictions
mutually  coincide.  

\begin{lem}\label{lem:37} We have
\begin{equation}\label{315}
\ltw{\sigma_n}(\pi \otimes\id)=(\id\otimes \pi )\tau
\end{equation}
for each $n\in\Bbb{Z}$.
\end{lem}

\begin{pf} Applying the appropriate twisting properties we obtain
\begin{equation*}
\begin{split}
\ltw{\sigma_n}(\pi \otimes \id)&=\ltw{\sigma_n}(\ilG\otimes \id)(\k\otimes
\id^2)(\phi\otimes \id)\\ 
&=(\id\otimes \ilG)(\sigma_n\otimes \id)(\id\otimes \sigma_n)(\k\otimes
\id^2)(\phi\otimes \id)\\ 
&=(\id\otimes\ilG)(\id\otimes\k\otimes\id)(\sigma_{-n}\otimes 
\id)(\id\otimes \sigma_n)(\phi \otimes  \id)\\
&=(\id\otimes  \ilG)(\id\otimes\k\otimes \id)(\id\otimes \phi)\tau= 
(\id\otimes \pi)\tau .\qed
\end{split}
\end{equation*}
\renewcommand{\qed}{}
\end{pf}

We are going to prove that the space $\Gamma$ is
naturally isomorphic to $\cal{A} \otimes \Gamma_{\inv}$, as a left
$\cal{A}$-module.

\begin{pro}\label{pro:38} Let us consider the map
$$
(\id\otimes P)\lac\colon\Gamma\rightarrow\cal{A} \otimes \Gamma_{\inv}.
$$
This is an isomorphism of left $\cal{A}$-modules. Its inverse is given by
$$
(\mlG\restr\cal{A} \otimes \Gamma_{\inv})\colon\cal{A} \otimes \Gamma_{\inv} 
\rightarrow\Gamma.
$$
\end{pro}

\begin{pf}
The map $\mlG\restr\cal{A}\otimes \Gamma_{\inv}$ is evidently a  left
$\cal{A}$-module homomorphism. Let us check that $\mlG\restr\cal{A} \otimes 
\Gamma_{\inv}$ and $(\id\otimes P)\lac$ are mutually  inverse  maps. 
Using \eqref{36}--\eqref{37} and \eqref{313} we obtain
\begin{multline*}
\mlG(\id\otimes P)\lac=\mlG (\id\otimes\mlG )(\id\otimes 
\k\otimes \id)(\id\otimes \lac)\lac\\ 
=\mlG(m\otimes \id)(\id\otimes \k\otimes \id)(\phi \otimes \id)\lac
=\mlG(\e1\otimes \id)\lac=\id.
\end{multline*}
On the other hand
$$P(a\vartheta)=\e(a)\vartheta$$
for each $a\in \cal{A}$  and $\vartheta\in \Gamma_{\inv}$. Using this and
\eqref{34} we obtain
\begin{multline*}
(\id\otimes P)\lac(a\vartheta)=(\id\otimes P)(m\otimes \mlG)(\id\otimes 
\sigma \otimes \id)(\phi (a)\otimes 1\otimes \vartheta)\\ 
=a^{(1)}\otimes P(a^{(2)}\vartheta)=a\otimes \vartheta.
\end{multline*}
Consequently,  the two maps are mutually inverse 
left $\cal{A}$-module isomorphisms.
\end{pf}

The  above  proposition allows us to identify $\Gamma\leftrightarrow
\cal{A} \otimes \Gamma_{\inv}$. In terms of this identification, the
following correspondences hold
\begin{gather}
d\leftrightarrow(\id\otimes \pi)\phi\label{317}\\
\lac\leftrightarrow\phi \otimes \id\label{316}\\
\mlG\leftrightarrow m\otimes \id.\label{318}
\end{gather}
The following technical lemma will be useful in some further computations. 

\begin{lem}\label{lem:39} We have
\begin{equation}\label{319}
 P[\pi(a)b]=\pi m\tau^{-1}\sigma (a\otimes b)-\e(a)\pi(b),
\end{equation}
for each $a,b\in\cal{A}$. 
\end{lem}

\begin{pf}
We compute
\begin{multline*}
P\mrG(\pi \otimes \id)=\mlG(\k\otimes \id)\lac
\mrG(\pi \otimes\id)\\
=\mlG(\k\otimes\id)(m\otimes\mrG)(\id\otimes
\ltw{\sigma}\otimes \id)(1\otimes \pi \otimes \phi).
\end{multline*}
According to \eqref{315}, this is further equal to
\begin{multline*}
\mlG(\k\otimes  \id)(m\otimes\mrG)(1\otimes\id\otimes   \pi   \otimes 
\id)(\tau \otimes \id)(\id\otimes \phi )=\\ 
=\mlG (\k\otimes\id)(\id\otimes\mrG)(\id\otimes\mlG \otimes 
\id)(\id\otimes \k\otimes d\otimes \id)(\id\otimes \phi \otimes 
\id)(\tau \otimes \id)(\id\otimes \phi )\\ 
=\mlG (\k\otimes \id)(\id\otimes \mlG)(\id\otimes \k\otimes dm)
(\id\otimes \phi\otimes \id)(\tau \otimes \id)(\id\otimes \phi )\\
-\mlG (\k\otimes \id)(\id\otimes \mlG)(\id^2\otimes\mlG)(\id\otimes
\k\otimes \id\otimes d)(\id\otimes \phi \otimes \id)(\tau \otimes  \id)
(\id\otimes\phi).
\end{multline*}
The first term in the above difference is transformed further
\begin{multline*}
\mlG(m\otimes \id)(\k\otimes\k\otimes  dm)(\id\otimes  \phi  \otimes 
\id)(\tau \otimes \id)(\id\otimes \phi)=\\
=\mlG (\k m\otimes  dm)(\tau^{-1}\sigma\tau^{-1}\sigma\tau^{-1} 
\otimes \id^2)(\id\otimes \phi  \otimes  \id)(\tau  \otimes 
\id)(\id\otimes \phi )\\ 
=\mlG (\k m\otimes  dm)(\tau^{-1}\otimes  \id^2)(\id\otimes  \phi  \otimes 
\id)(\sigma \tau^{-1}\sigma \otimes \id)(\id\otimes \phi )\\
=\mlG(\k m\otimes dm)(\id\otimes \sigma \otimes \id)(\phi \otimes  \phi 
)\tau^{-1}\sigma =\mlG (\k\otimes  d)\phi m\tau^{-1}\sigma=\pi m\tau^{-1}
\sigma .
\end{multline*}
Concerning the second term,
$$
\mlG (\k\otimes \id)(\id\otimes \mlG )(\id\otimes 1\e\otimes
d)(\tau  \otimes\id)(\id\otimes \phi )=\mlG (\k\otimes d)(\e\otimes
\phi )=\e\otimes  \pi,
$$
which completes the proof.
\end{pf}

Let $\cal{R}$ be the intersection of spaces $\ker(\pi)$ and $\ker(\e)$.
As follows directly from the previous lemma, the space $\cal{R}$ is a 
right ideal in the algebra $\cal{A}_0$,  which  coincides  as  a 
vector space with $\cal{A}$, but which is endowed with the  product 
$m_0=m\tau^{-1}\sigma$, as discussed in \cite{D-BQG}-Appendix.
According to \eqref{315}, we have
\begin{equation}\label{320}
 \tau (\cal{R} \otimes \cal{A} )=\cal{A} \otimes \cal{R} .
\end{equation}
The map $\pi$ induces the isomorphism
\begin{equation}
\Gamma_{\inv}\leftrightarrow\ker(\e)/\cal{R}.
\end{equation}

It is easy to see that the map ${\circ}\colon\Gamma_{\inv}\otimes \cal{A}
\rightarrow\Gamma_{\inv}$ given by
\begin{equation}\label{321}
\pi(a)\circ b=P(\pi(a)b)=\pi m_0(a\otimes b)
\end{equation}
defines a right $\cal{A}_0$-module structure on the space $\Gamma_{\inv}$. 
In the above formula it is assumed that $a\in \ker(\e)$,
while $b$ is arbitrary.

In terms of the identification $\Gamma\leftrightarrow\cal{A}\otimes
\Gamma_{\inv}$ the right $\cal{A}$-module structure is given by
\begin{equation*}\label{322}
\mrG\leftrightarrow(m\otimes{\circ})(\id\otimes  \sigma_*\otimes\id)
(\id^2\otimes\phi),
\end{equation*}
where $\sigma_*\colon\Gamma_{\inv}\otimes \cal{A}\rightarrow\cal{A} \otimes
\Gamma_{\inv}$ is the common  left-invariant  part of all operators
$\ltw{\sigma_n}$. We shall now prove that $\Gamma$ is trivial as
a right $\cal{A}$-module.

\begin{pro}\label{pro:310} The multiplication map
\begin{equation}\label{323}
(\mrG\restr\Gamma_{\inv}\otimes  \cal{A})\colon\Gamma_{\inv}\otimes\cal{A} 
\rightarrow\Gamma=\cal{A} \otimes \Gamma_{\inv}
\end{equation}
is an isomorphism of right $\cal{A}$-modules. Its inverse is given by
\begin{equation}\label{324}
 ({\circ}\otimes\k)(\id\otimes  \phi \k^{-1})\sigma^{-1}_*\colon
\cal{A}\otimes\Gamma_{\inv}\rightarrow\Gamma_{\inv}\otimes\cal{A}.
\end{equation}
\end{pro}

\begin{pf}
Clearly, \eqref{323} is a right $\cal{A}$-module homomorphism. A
direct computation gives
\begin{multline*}
\bigl[(\id\otimes {\circ})(\sigma_*\otimes \id)({\circ}\otimes
\phi\k)(\id\otimes\phi\k^{-1})\sigma^{-1}_*\bigr](\id\otimes \pi)=\\ 
=(\id\otimes{\circ})(\sigma_*\otimes \id)({\circ}\otimes \phi\k)(
\pi  \otimes\phi\k^{-1})\tau^{-1}\\ 
=(\id\otimes {\circ})(\sigma_*\otimes  \id)(\pi  m\tau^{-1}\sigma
\otimes  \phi\k)(\id\otimes \phi \k^{-1}  )\tau^{-1}\\ 
=(\id\otimes  \pi)(\id\otimes   m\tau^{-1}\sigma)(\tau\otimes 
\id)(m\otimes \phi\k)(\id\otimes  \phi)\tau^{-1}\sigma(\id\otimes 
\k^{-1})\tau^{-1}\\
=(\id\otimes   \pi)(\id\otimes    m\tau^{-1})(\phi\otimes\id)\sigma
(m\otimes\k)(\id\otimes   \phi   \k^{-1}  )\sigma^{-1}\\ 
=(\id\otimes  \pi  )(\id\otimes   m)(\sigma 
\otimes \id)\bigl[m\otimes  \sigma  (\k\otimes 
\k)\phi\bigr](\id\otimes \phi \k^{-1})\sigma^{-1}\\
=(\id\otimes  \pi)\sigma(m\otimes   \id)(m\otimes\k\otimes 
\k)(\id\otimes \phi \otimes \id)(\id\otimes \phi \k^{-1})\sigma^{-1}\\
=(\id\otimes  \pi)\sigma(m\otimes   \id)(\id\otimes   1\e\otimes 
\k)(\id\otimes \phi \k^{-1})\sigma^{-1}=(\id\otimes \pi).
\end{multline*}
Furthermore,
\begin{multline*}
({\circ}\otimes \k)(\id\otimes \phi  \k^{-1})\sigma_*^{-1}
(\id\otimes {\circ})(\sigma_*\otimes \id)(\pi \otimes \phi)=\\
=({\circ}\otimes \k)(\id\otimes  \phi  \k^{-1}  )\sigma_*^{-1}
(\id\otimes  \pi m\tau^{-1}\sigma)(\tau \otimes \id)(\id\otimes \phi)\\
=(\pi m\tau^{-1}\sigma   \otimes\k)(\id\otimes   \phi   \k^{-1})\tau^{-1}
(\id\otimes m\tau^{-1}\sigma )(\tau \otimes \id)(\id\otimes \phi )\\
=(\pi m\otimes \k)(\id\otimes \phi \k^{-1}  )\sigma^{-1}
(\id\otimes  m)(\sigma\otimes \id)(\id\otimes \phi )\\
=(\pi m\otimes \k)(m\otimes \phi \k^{-1})(\id\otimes \sigma^{-1})
(\id\otimes \phi )\\
=(\pi m\otimes\id)(\id\otimes m\otimes\id)(\id\otimes\k\otimes 
\id\otimes\k)(\id\otimes \phi \otimes \id)(\id\otimes \phi \k^{-1})\\
=(\pi m\otimes \id)(\id\otimes 1\e\otimes \k)(\id\otimes \phi  \k^{-1})=
(\pi\otimes \id).
\end{multline*}
The above computations  are performed in  the  spaces
$\cal{A}\otimes\ker(\e)$ and $\ker(\e)\otimes \cal{A}$ respectively.
\end{pf}

 In the framework  of the identification  $\Gamma\leftrightarrow
\Gamma_{\inv}\otimes\cal{A}$, the following correspondences hold:
\begin{gather}
-d \leftrightarrow (\pi \k^{-1}\otimes \id)\sigma^{-1}\phi =(\pi \otimes
\k)\phi \k^{-1}\label{328}\\
\mrG \leftrightarrow\id\otimes m \label{325}\\ 
 \lac \leftrightarrow (\sigma_*\otimes \id)(\id\otimes \phi)\label{326}\\ 
\mlG\leftrightarrow\bigl[{\circ}(\id\otimes  \k^{-1})\otimes  m\bigr]
(\id\otimes  \sigma^{-1}\phi
\otimes \id)(\sigma_*^{-1}\otimes \id).\label{327}
\end{gather}
These correspondences follow from  \eqref{323}--\eqref{324}, performing
simple algebraic transformations.

We are ready to  present  a  braided  counterpart  of  the
 reconstruction theorem \cite{W} of the standard theory.  As  we  have 
seen, every left-covariant  calculus  $\Gamma$  is  completely 
determined by the corresponding $\cal{R}$.
The  following  proposition  shows  that conversely,
every right $\cal{A}_0$-ideal which  satisfies \eqref{320}
naturally gives rise to a first-order left-covariant calculus.

\begin{pro}\label{pro:311} Let $\cal{R}\subseteq\ker(\e)$  be
an arbitrary $\tau$-invariant right $\cal{A}_0$-ideal. 
Let us define spaces $\Gamma_{\inv}$ and $\Gamma$,
together with maps ${\circ}\colon\Gamma_{\inv}\otimes\cal{A}
\rightarrow\Gamma_{\inv}$ and $\pi\colon\cal{A}\rightarrow
\Gamma_{\inv}$, as well as $\sigma_*\colon\Gamma_{\inv}
\otimes \cal{A}\rightarrow\cal{A} \otimes \Gamma_{\inv}$ by the equalities
\begin{gather*}
\Gamma_{\inv}=\ker(\e)/\cal{R}\qquad\Gamma=\cal{A}\otimes \Gamma_{\inv}\\     
\pi (a)=\bigl[a-\e(a)\bigr]_{\cal{R}}\\           
(\pi (a)\circ b)=\pi\bigl[m_0 (a\otimes b)-\e(a)b\bigr]\\  
\sigma_*(\pi \otimes \id)=(\id\otimes \pi )\tau.
\end{gather*}
Finally, let us define maps
$$\lac\colon\Gamma \rightarrow\cal{A} \otimes  \Gamma\qquad
d\colon\cal{A}\rightarrow\Gamma\qquad\mlG\colon\cal{A}  \otimes  \Gamma
\rightarrow\Gamma\qquad
\mrG\colon\Gamma  \otimes  \cal{A}\rightarrow\Gamma$$
by equalities \eqref{316}--\eqref{318}  and \eqref{322}
respectively.

 Then,  $\mlG$ and  $\mrG$   determine  a  structure  of  a
unital $\cal{A}$-bimodule on $\Gamma$. Moreover, $\Gamma$ is a
left-covariant first-order differential calculus over $G$,
with the differential and the left action coinciding with
the introduced $d$ and $\lac$ respectively.
\end{pro}

\begin{pf}
It is clear that $\mlG$ determines a left
$\cal{A}$-module structure on $\Gamma$. Let us prove that $\mrG$
determines a right $\cal{A}$-module  structure. We have
\begin{multline*}
\mrG(\mrG\otimes\id)=(m\otimes{\circ})(\id\otimes\sigma_*\otimes 
\id)(\id^2 \otimes  \phi)\bigl[(m\otimes{\circ})(\id\otimes
\sigma_*\otimes\id)(\id^2\otimes \phi )\otimes \id\bigr]\\ 
=(m\otimes {\circ})(\id^2 \otimes {\circ}\otimes  \id)(\id\otimes  \sigma_* 
\otimes\id^2)(m\otimes\id\otimes\tau\otimes 
\id)(\id\otimes\sigma_*\otimes\id\otimes\phi 
)(\id^2 \otimes \phi \otimes \id)\\
=(m\otimes{\circ})(\id\otimes   \sigma_*\otimes   m\tau^{-1}\sigma 
)(m\otimes \id\otimes \tau \otimes \id)(\id\otimes \sigma_*\otimes 
\id\otimes \phi )(\id^2\otimes \phi \otimes \id)\\
=(m\otimes {\circ})(m\otimes  \sigma_*\otimes  m)(\id^3 \otimes  \sigma 
\otimes   \id)(\id\otimes   \sigma_*\otimes    \id\otimes    \phi 
)(\id^2 \otimes \phi \otimes \id)\\
=(m\otimes {\circ})(\id\otimes m\otimes \id^2)(\id^2\otimes
\sigma_*\otimes m)(\id\otimes \sigma_*\otimes \sigma \otimes \id)(\id^2
\otimes  \phi\otimes \phi )\\
=(m\otimes {\circ})(\id\otimes \sigma_*\otimes \id)(\id^2\otimes  m\otimes 
m)(\id^3 \otimes \sigma \otimes \id)(\id^2 \otimes \phi \otimes \phi )\\
=(m\otimes{\circ})(\id\otimes  \sigma_*\otimes  \id)(\id^2 \otimes
\phi m)= \mrG(\id\otimes m).
\end{multline*}

Here,  we  have  used \eqref{321}--\eqref{322}, and identities
\begin{gather}
\sigma_*({\circ}\otimes   \id)=(\id\otimes{\circ})(\sigma_*\otimes 
\id)(\id\otimes \tau )\label{332}\\ 
(m\otimes \id)(\id\otimes \sigma_*)(\sigma_*\otimes  \id)=\sigma_*
(\id\otimes m),\label{333}
\end{gather}
which follow from \eqref{315} and \eqref{321}.

    The maps $\mlG$ and $\mrG$ commute, because
\begin{multline*}
\mrG(\mlG\otimes\id)= \mrG(m\otimes 
\id^2)=(m\otimes \id)(m\otimes \id\otimes {\circ})(\id^2\otimes
\sigma_*\otimes\id)(\id^3 \otimes \phi)\\ 
=(m\otimes \id)\bigl[\id\otimes (m\otimes {\circ})(\id\otimes
\sigma_*\otimes\id)(\id^2 \otimes \phi)\bigr]=\mlG(\id\otimes\mrG).
\end{multline*}
It is easy to see that the bimodule $\Gamma$ is unital.

    According to Lemma~\ref{lem:39} and equation~\eqref{321},
\begin{equation}\label{334}
 \pi m=[\e\otimes \pi +{\circ}(\pi \otimes \id)]\sigma^{-1}\tau.
\end{equation}
Using this, equations \eqref{317} and \eqref{318} and \eqref{322}
we obtain
\begin{multline*}
dm=(m\otimes \pi m)(\id\otimes  \sigma  \otimes  \id)(\phi  \otimes 
\phi )=\\=\bigl[m\otimes (\e\otimes \pi )\sigma^{-1}\tau\bigr]
(\id\otimes  \sigma\otimes \id)(\phi \otimes \phi )+ 
\bigl[m\otimes {\circ}(\pi  \otimes  \id)\sigma^{-1}\tau\bigr]
(\id\otimes  \sigma\otimes \id)(\phi \otimes \phi )\\
=(m\otimes  \e\otimes  \pi  )(\id\otimes  \tau  \otimes  \id)(\phi 
\otimes \phi  )+\bigl[m\otimes {\circ}(\pi  \otimes  \id)\bigr]
(\id\otimes  \tau\otimes \id)(\phi \otimes \phi )\\
=(m\otimes\pi )(\id\otimes\phi)+(m\otimes 
{\circ})(\id\otimes \sigma_*\otimes \id)(\id\otimes \pi \otimes \id^2)
(\phi \otimes \phi)\\=\mlG (\id\otimes d)+ \mrG(d\otimes \id).
\end{multline*}

To complete the proof, let us observe that \eqref{317} implies that
$\lac$ given by \eqref{316} is indeed the left action.
\end{pf}

\section{Bicovariant Calculi}

 In this section we shall study bicovariant differential calculi
$\Gamma$ over $G$. As in the standard theory \cite{W}
the right action  $\rac$ and the left action $\lac$ are mutually
compatible.

\begin{pro}\label{pro:41} The diagram
\begin{equation}\label{41}
\begin{CD}
\Gamma @>{\mbox{$\lac$}}>> \cal{A}\otimes\Gamma\\
@V{\mbox{$\rac$}}VV @VV{\mbox{$\id\otimes\rac$}}V\\
\Gamma\otimes\cal{A} @>>{\mbox{$\lac\otimes\id$}}> \cal{A}\otimes\Gamma
\otimes\cal{A}
\end{CD}
\end{equation}
is commutative. 
\end{pro}

\begin{pf} Applying \eqref{32} and \eqref{A1} we obtain
\begin{equation*}
\begin{split}
(\lac\otimes \id)& \rac\ilG=(\lac   \otimes  \id)(\ilG\otimes 
m)(\id\otimes \sigma \otimes \id)(\phi \otimes \phi )\\ 
&=(m\otimes \ilG\otimes m)(\id\otimes  \sigma  \otimes  \id^3)(\phi 
\otimes  \phi  \otimes  \id^2)(\id\otimes  \sigma  \otimes  \id)(\phi 
\otimes \phi )\\
&=(m\otimes    \ilG\otimes    m)(\id^3 \otimes    \sigma     \otimes 
\id)(\id^2 \otimes  \phi  \otimes  \phi  )(\id\otimes  \sigma   \otimes 
\id)(\phi \otimes \phi )\\
&=(m\otimes  \rac\ilG )(\id\otimes  \sigma  \otimes  \id)(\phi 
\otimes \phi )=(\id\otimes  \rac)\lac\ilG. \qed
\end{split}
\end{equation*}
\renewcommand{\qed}{}
\end{pf}

As a simple consequence of \eqref{41} we find that
the spaces $\Gamma_{\inv}$ and $\Ginvr$ are right/left-invariant
respectively. 
The  following proposition  characterizes   the corresponding
restrictions of  $\rac$ and $\lac$. 
Let $\ad\colon\cal{A}\rightarrow\cal{A}\otimes\cal{A}$  be  the  adjoint 
action of $G$ on itself, as defined in Appendix B.

\begin{pro}\label{pro:42} The following identities hold
\begin{align}
\rac\pi&=(\pi \otimes \id)\ad\label{43a}\\                             
\lac\pir  &=(\id\otimes \pir)\tau (\k\otimes \k)\ad\k^{-1}.\label{43b}
\end{align}
\end{pro}

\begin{pf} We compute
\begin{equation*}
\begin{split}
 \rac\pi&= \rac\ilG(\k\otimes   \id)\phi=(\ilG\otimes 
m)(\id\otimes  \sigma  \otimes  \id)(\sigma  \otimes   \id^2)(\k\otimes 
\k\otimes \id^2)(\phi \otimes \phi )\phi\\ 
&=(\ilG\otimes m)(\k\otimes \id\otimes \k\otimes \id)(\id\otimes  \tau 
\sigma^{-1}\tau \otimes \id)(\sigma \otimes  \id^2)(\phi  \otimes\phi)\phi\\
&=(\ilG\otimes m)(\k\otimes  \id\otimes  \k\otimes  \id)(\phi  \otimes 
\id^2)(\tau\otimes   \id)(\id\otimes   \phi)\phi\\
&=(\pi\otimes m)(\id\otimes \k\otimes \id)
(\tau \otimes \id)(\id\otimes \phi)\phi=(\pi \otimes \id)\ad.
\end{split}
\end{equation*}
Completely similarly,
\begin{equation*}
\begin{split}
 \lac\pir    &=\lac\irG(\id\otimes\k)\phi=(m\otimes 
\irG)(\id\otimes    \sigma    \otimes\id)(\id^2 \otimes\sigma)(\id^2
\otimes \k\otimes \k)(\phi \otimes \phi )\phi\\
&=(m\otimes\irG)(\id\otimes\k\otimes    \id\otimes 
\k)(\id\otimes \tau \sigma^{-1}\tau  \otimes  \id)(\id^2 \otimes 
\sigma )(\phi \otimes \phi )\phi\\ 
&=(m\otimes  \irG)(\id\otimes \k\otimes \id\otimes  \k)(\id^2 \otimes 
\phi )(\id\otimes \tau )(\phi \otimes \id)\phi\\
&=(m\otimes \pir)(\id\otimes \k\otimes \id)(\id\otimes \tau )(\phi  \otimes 
\id)\phi =(\id\otimes \pir)\tau(\k\otimes \k)\ad\k^{-1}.\qed
\end{split}
\end{equation*}
\renewcommand{\qed}{}
\end{pf}

We pass to the the analysis of the specific
twisting properties of the left and the right action maps.

\begin{pro}\label{pro:43} The following equalities hold
\begin{align}
(\lac\otimes\id)\rtw{\sigma_{n+m}}&=
(\id\otimes\rtw{\sigma_m})(\sigma_n\otimes\id)(\id\otimes\lac)\label{44a}\\
(\id\otimes\rac)\ltw{\sigma_{n+m}}&=(\ltw{\sigma_m}\otimes\id)(\id\otimes
\sigma_n)(\rac\otimes\id)\label{44b}
\end{align}
for each $n,m\in\Bbb{Z}$. 
\end{pro}

\begin{pf} A direct computation gives
\begin{multline*}
(\id\otimes\rtw{\sigma_m})(\sigma_n\otimes   \id)(\id\otimes \lac)
(\id\otimes \ilG)=\\
=(\id\otimes\rtw{\sigma_m})(\sigma_n\otimes 
\id)(\id\otimes m\otimes \ilG)(\id^2\otimes \sigma 
\otimes \id)(\id\otimes \phi \otimes \phi )\\ 
=(m\otimes\rtw{\sigma_m})(\id^3\otimes \ilG)(\id\otimes \sigma_n
\otimes \id^2 )(\sigma_n  \otimes
\id^3)(\id^2\otimes \sigma \otimes \id)(\id\otimes \phi \otimes \phi )\\
=(m\otimes \ilG\otimes\id)(\id^3 \otimes \sigma_m)(\id^2 \otimes \sigma_m
\otimes \id)(\id\otimes \sigma_n  \otimes
\id^2)(\sigma_n  \otimes \sigma \otimes \id)(\id\otimes \phi \otimes \phi )\\
=(m\otimes \ilG\otimes\id)(\id\otimes \sigma \otimes\sigma_m)(\id^2 \otimes 
\sigma_n  \otimes \id)(\id\otimes \sigma_m\otimes \id^2)(\sigma_n\otimes\id^3)
(\id\otimes \phi\otimes \phi )\\ 
=(m\otimes \ilG\otimes \id)(\id\otimes  \sigma  \otimes  \id^2)(\phi 
\otimes \phi \otimes \id)(\id\otimes \sigma_{n+m})(\sigma_{n+m}\otimes\id)\\ 
=(\lac\ilG\otimes  \id)(\id\otimes  \sigma_{n+m})(\sigma_{n+m}\otimes\id)
=(\lac\otimes\id) \rtw{\sigma_{n+m}}(\id\otimes\ilG).
\end{multline*}
Similarly we obtain
\begin{multline*}
(\ltw{\sigma_m}\otimes \id)(\id\otimes  \sigma_n)(\rac  \otimes 
\id)( \irG\otimes \id)=\\
=(\ltw{\sigma_m}\otimes \id)(\id\otimes \sigma_n)( \irG\otimes  m\otimes 
\id)(\id\otimes \sigma \otimes \id^2)(\phi \otimes \phi \otimes \id)\\
=(\ltw{\sigma_m}\otimes m)(\id\otimes \sigma_n\otimes  \id)( \irG\otimes 
\id\otimes \sigma_n)(\id\otimes  \sigma  \otimes  \id^2)(\phi  \otimes 
\phi \otimes \id)\\
=(\id\otimes\irG\otimes m)(\sigma_m  \otimes \id^3)(\id\otimes
\sigma_m \otimes\id^2)(\id^2\otimes  \sigma_n\otimes\id)(\id\otimes
\sigma   \otimes\sigma_n)(\phi  \otimes \phi \otimes \id)\\
=(\id\otimes \irG\otimes m)(\sigma_m  \otimes \sigma \otimes \id)(\id\otimes 
\sigma_n      \otimes\id^2)(\id^2 \otimes      \sigma_m\otimes 
\id)(\id^3\otimes   \sigma_n)(\phi \otimes \phi \otimes \id)\\ 
=(\id\otimes\irG\otimes   m)(\sigma_m\otimes\sigma    \otimes 
\id)(\id\otimes \sigma_n\otimes \phi )(\phi \otimes \sigma_{n+m})\\
=(\id\otimes  \irG\otimes m)(\id^2\otimes \sigma \otimes \id)(\id\otimes 
\phi \otimes \phi )(\sigma_{n+m}\otimes \id)(\id\otimes \sigma_{n+m})\\
=(\id\otimes\rac\irG)(\sigma_{n+m}\otimes \id)(\id\otimes \sigma_{n+m})= 
(\id\otimes\rac)\ltw{\sigma_{n+m}}(\irG\otimes \id).\qed
\end{multline*}
\renewcommand{\qed}{}
\end{pf}

As a simple consequence of the previous proposition we find
\begin{align}
\rtw{\sigma_n}(\cal{A} \otimes \Gamma_{\inv})&=\Gamma_{\inv}\otimes
\cal{A}\label{45a}\\
\ltw{\sigma_n}(\Ginvr\otimes \cal{A} )&=\cal{A} \otimes\Ginvr.\label{45b}
\end{align}
The  following  proposition  describes  the corresponding
restriction twistings. 

\begin{pro}\label{pro:45} The following identities hold
\begin{gather}
\rtw{\sigma_n}(\id\otimes \pi )=(\pi \otimes \id)\tau \label{46a}\\
\ltw{\sigma_n}(\pir\otimes \id)=(\id\otimes \pir)\tau .\label{46b}
\end{gather}
\end{pro}

\begin{pf}
Using standard twisting transformations we obtain
\begin{equation*}
\begin{split}
\rtw{\sigma_n}(\id\otimes\pi )&=(\ilG\otimes 
\id)(\id\otimes \sigma_n)(\sigma_n  \otimes 
\id)(\id\otimes\k\otimes \id)(\id\otimes \phi)\\ 
&=(\ilG\otimes\id)(\k\otimes\id^2 )(\id\otimes        \sigma_n)
(\sigma_{-n}\otimes \id)(\id\otimes \phi)\\
&=(\ilG\otimes \id)(\k\otimes \id^2)(\phi \otimes \id)\tau=(\pi \otimes
\id)\tau.
\end{split}
\end{equation*}
The second identity can be derived in a similar manner.
\end{pf}

Let $\cal{R}\subseteq\ker(\e)$ be the right $\cal{A}_0$-ideal which canonically
corresponds to $\Gamma$. In
the following proposition we have characterized
bicovariance in terms of $\cal{R}$.

\begin{pro}\label{pro:46} \bla{i} We have
\begin{gather}
\ad(\cal{R} )\subseteq\cal{R} \otimes \cal{A} \label{47}\\ 
\tau (\cal{A} \otimes \cal{R} )=\cal{R} \otimes \cal{A} . \label{48}
\end{gather}

\bla{ii} Conversely, if $\cal{R}\subseteq\ker(\e)$
corresponding to a left-covariant calculus $\Gamma$ is $\ad$-invariant, 
then the calculus $\Gamma$ is bicovariant. Moreover,  in  terms
of the identification $\Gamma =\Gamma_{\inv}\otimes \cal{A}$,  
the right action  $\rac\colon\Gamma\rightarrow\Gamma \otimes \cal{A}$
is given by
\begin{equation}\label{49}
\rac=(\id^2 \otimes m)(\id\otimes \sigma \otimes \id)(\adj\otimes \phi),
\end{equation}
where the map $\adj\colon\Gamma_{\inv}\rightarrow
\Gamma_{\inv}\otimes\cal{A}$ is given by
\begin{equation}\label{410}
                     \adj\pi =(\pi \otimes \id)\ad.                           
\end{equation}
\end{pro}

\begin{pf} The  first  statement  of  the  proposition  is  a  direct 
consequence of \eqref{43a} and \eqref{46a}.
Concerning the second part, it is
sufficient to check that the map $\xi$ given by the right-hand side of
\eqref{49} satisfies \eqref{A2}--\eqref{A3}.
Using  the  structuralization $\Gamma=\Gamma_{\inv}\otimes\cal{A}$
as well as equalities \eqref{328} and \eqref{410} we obtain
\begin{equation*}
\begin{split}
 \xi d&=-(\id^2 \otimes m)(\id\otimes\sigma\otimes 
\id)(\adj\otimes \phi )(\pi \otimes \id)(\id\otimes \k)\phi k^{-1}\\
&=-(\pi \otimes \id\otimes m)(\id\otimes \sigma \otimes  \id)\bigl[\ad\otimes 
\sigma (\k\otimes \k)\phi\bigr]\phi \k^{-1}\\
&=-(\pi  \otimes  \sigma)(\id\otimes   m\otimes\id)\bigl(\ad\otimes 
(\k\otimes \k)\phi \bigr)\phi k^{-1}\\
&=-(\pi   \otimes\sigma)(\id\otimes   \k\otimes   \k)(\tau   \otimes 
\id)(\id\otimes \phi )\phi k^{-1}\\
&=-(\pi \otimes \k\otimes  \k)(\phi  \otimes\id)\sigma \phi k^{-1}\\ 
&=-(\pi  \otimes  \k\otimes\id)(\phi\k^{-1}  \otimes  \id)\phi
=(d\otimes\id)\phi.
\end{split}
\end{equation*}
Furthermore, \eqref{325} implies
\begin{multline*}
\xi\mrG=(\id^2\otimes m)(\id\otimes \sigma \otimes \id)(\adj\otimes\phi m)=\\ 
=(\id^2\otimes m)(\id\otimes\sigma\otimes\id)(\adj\otimes m\otimes m)(\id^2
\otimes\sigma\otimes \id)(\id\otimes \phi \otimes \phi )\\ 
=(\id\otimes m\otimes m)(\id^2\otimes\sigma \otimes 
\id)(\id\otimes\sigma \otimes\id^2)(\adj\otimes\id^2\otimes 
m)(\id^2 \otimes \sigma \otimes \id)(\id\otimes \phi \otimes \phi )\\ 
=(\id\otimes m\otimes m)(\id^3 \otimes 
m\otimes \id)(\id^2 \otimes \sigma \otimes \id^2 )(\id\otimes \sigma \otimes 
\sigma \otimes \id)(\adj\otimes \phi \otimes \phi )\\ 
=(\id\otimes m\otimes m)(\id^2 \otimes  \sigma  \otimes  \id)(\id^2 \otimes 
m\otimes \id^2)(\id\otimes \sigma \otimes \id^3)(\adj\otimes \phi \otimes 
\phi )\\ 
=(\mrG\otimes m)(\id\otimes \sigma \otimes \id)(\xi\otimes \phi ).
\end{multline*}
Consequently, $\Gamma$  is   bicovariant   and $\xi=\rac$.
\end{pf}

\section{Antipodally Covariant Calculi}

    In this Section we  shall consider differential  structures 
covariant relative to the antipode map.

\begin{defn}\label{def:51} A  first-order  calculus $\Gamma$
is  called $\k$-covariant iff the following equivalence holds
$$
 \omega \in \ker(\ilG) \iff \omega \in \ker\bigl[\irG(\k\otimes  \k)\tau 
\sigma^{-1}\tau \sigma^{-1}\tau\bigr].
$$
\end{defn}

    Let us assume that $\Gamma$ is $\k$-covariant. Then the formula
\begin{equation}\label{51}
 \K\ilG = \irG(\k\otimes \k)\tau\sigma^{-1}\tau\sigma^{-1}\tau
\end{equation}
consistently and  uniquely  determines  a  bijective  map
$\K\colon\Gamma\rightarrow\Gamma$. It follows that
\begin{align}
                 d\k&=\K d\label{52}\\
\K\irG&=\ilG(\k\otimes\k)\tau\sigma^{-1}\tau\sigma^{-1}\tau .\label{53}
\end{align}
Let us analyze properties of $\Gamma$,  in  the  case
when it is also $\sigma$-covariant. 

\begin{pro}\label{pro:51}
\bla{i} If $\Gamma$ is left $\sigma$-covariant  (and accordingly, 
left $\cal{F}$-covariant) then
\begin{gather}
\ltw{\sigma_n}(\K\otimes \id)=(\id\otimes \K)\ltw{\sigma_{-n}}\label{55} \\
\K\mrG=\mlG (\k\otimes \K)\ltw{(\tau\sigma^{-1}\tau\sigma^{-1}\tau)}.\label{56}
\end{gather}

\smallskip
\bla{ii} If $\Gamma$ is right $\cal{F}$-covariant then
\begin{gather}
\rtw{\sigma_n}(\id\otimes \K)=(\K\otimes \id)\rtw{\sigma_{-n}}\label{57}\\
\K\mlG=\mrG(\K\otimes \k)\rtw{(\tau \sigma^{-1}\tau\sigma^{-1}\tau)}.\label{58} 
\end{gather}
\end{pro}

\begin{pf} Let us assume that $\Gamma$ is left $\cal{F}$-covariant. A direct
computation gives
\begin{equation*}
\begin{split}
\ltw{\sigma_n}(\K\otimes\id)&(\ilG\otimes\id)=\ltw{\sigma_n}(\irG\otimes 
\id)(\k\otimes \k\otimes \id)(\tau \sigma^{-1}\tau \sigma^{-1}\tau
\otimes \id)\\
&=(\id\otimes  \irG)(\sigma_n\otimes \id)(\id\otimes \sigma_n)(\k\otimes 
\k\otimes \id)(\tau \sigma^{-1}\tau \sigma^{-1}\tau \otimes \id)\\ 
&=(\id\otimes\irG)(\id\otimes\k\otimes 
\k)(\id\otimes \tau \sigma^{-1}\tau \sigma^{-1}\tau)(\sigma_{-n}
\otimes \id)(\id\otimes \sigma_{-n})\\
&=(\id\otimes \K\ilG)(\sigma_{-n}\otimes 
\id)(\id\otimes \sigma_{-n})=(\id\otimes  \K)\ltw{\sigma_{-n}}
(\ilG\otimes \id).
\end{split}
\end{equation*}
Furthermore,
\begin{equation*}
\begin{split}
\K\mrG(\irG\otimes \id)&=\ilG(\k\otimes\k)\sigma_{-2}(\id\otimes 
m)=\ilG(\k m\otimes \k)(\id\otimes \sigma_{-2})(\sigma_{-2}\otimes \id)\\ 
&=\ilG(m\otimes   \id)(\k\otimes   \k\otimes   \k)(\sigma_{-2}\otimes 
\id)(\id\otimes \sigma_{-2})(\sigma_{-2}\otimes \id)\\ 
&=\mlG (\id\otimes   \ilG)(\k\otimes\k\otimes\k)(\id\otimes\sigma_{-2} 
)(\sigma_{-2}\otimes \id)(\sigma_{-2}\otimes \id)\\
&=\mlG (\k\otimes  \K \irG)(\sigma_{-2}\otimes   \id)(\sigma_{-2}\otimes 
\id)=\mlG(\k\otimes \K)\ltw{\sigma_{-2}}(\ilG\otimes \id).
\end{split}
\end{equation*}
Symmetrically, assuming the right $\cal{F}$-covariance of $\Gamma$ we
get
\begin{equation*}
\begin{split}
\rtw{\sigma_n}(\id\otimes \K)(\id\otimes
\irG)&=\rtw{\sigma_n}(\id\otimes\ilG)(\id\otimes\k\otimes 
\k)(\id\otimes \sigma_{-2})\\ 
&=(\ilG\otimes    \id)(\id\otimes    \sigma_n)(\sigma_n\otimes 
\id)(\id\otimes \k\otimes \k)(\id\otimes \sigma_{-2})\\
&=(\ilG\otimes   \id)(\k\otimes   \k\otimes   \id)(\sigma_{-2}\otimes 
\id)(\id\otimes \sigma_{-n})(\sigma_{-n}\otimes \id)\\ 
&=(\K\otimes  \id)(\irG\otimes   \id)(\id\otimes   \sigma_{-n})(\sigma_{-n} 
\otimes \id)\\&=(\K\otimes \id)\rtw{\sigma_{-n}}(\id\otimes\irG).
\end{split}
\end{equation*}
Finally,
\begin{equation*}
\begin{split}
   \K\mlG(\id\otimes \ilG)&= \irG(\k\otimes\k)\sigma_{-2}(m\otimes\id)=
\irG(\k\otimes\k m)(\sigma_{-2}\otimes \id)(\id\otimes \sigma_{-2})\\
&=\irG(\id\otimes m)(\k\otimes\k\otimes\k)(\id\otimes\sigma_{-2} 
)(\sigma_{-2}\otimes \id)(\id\otimes \sigma_{-2})\\
&=\mrG(\irG\otimes   \id)(\k\otimes   \k\otimes   \k)(\sigma_{-2}\otimes 
\id)(\id\otimes \sigma_{-2})(\sigma_{-2}\otimes\id)\\ 
&=\mrG(\K\ilG\otimes\k)(\id\otimes\sigma_{-2})(\sigma_{-2}\otimes 
\id)= \mrG(\K\otimes \k) \rtw{\sigma_{-2}}(\id\otimes \ilG),
\end{split}
\end{equation*}
which completes the proof. 
\end{pf}

Now, we shall analyze interrelations  between  $\k$-covariance  and
bicovariance.

\begin{pro}\label{pro:52} A left-covariant calculus $\Gamma$ 
is $\k$-covariant if and only if it is bicovariant. 
In this case the following identities hold:
\begin{align}
 \lac \K&=(\k\otimes \K)\ltw{\sigma}\rac \label{59}\\
  \rac \K&=(\K\otimes \k)\rtw{\sigma}\lac  . \label{510}
\end{align}
Moreover, the diagram 
\begin{equation}\label{511}
\begin{CD}
\cal{A}\otimes\Gamma\otimes\cal{A} @>{\mbox{$\k\otimes\id\otimes\k$}}>>
\cal{A}\otimes\Gamma\otimes\cal{A} \\
@AAA @VVV\\
\Gamma @>>{\mbox{$-\K$}}> \Gamma
\end{CD}
\end{equation}
is commutative. Here, the vertical arrows are the corresponding
double-sided actions and products. 
\end{pro}

\begin{pf}
Let  us assume  that $\Gamma$ is left-covariant and 
$\k$-covariant, and  let us consider  a  map
$\xi\colon\Gamma\rightarrow\Gamma\otimes  \cal{A}$
defined by
$$
\xi=(\K^{-1}\otimes \k^{-1})(\ltw{\sigma})^{-1}\lac \K.
$$
It turns out that $\xi$ is the right action for $\Gamma$. Indeed, 
\begin{multline*}
\xi\irG=(\K^{-1}\otimes \k^{-1})(\ltw{\sigma })^{-1}\lac\ilG (\k\otimes \k)
\sigma_{-2}=\\
=(\K^{-1}\otimes\k^{-1})(\ltw{\sigma})^{-1}(m\otimes\ilG)(\id\otimes   \sigma 
\otimes \id)(\phi \k\otimes \phi \k)\sigma_{-2}\\
=(\K^{-1}\otimes\k^{-1})(\ilG \otimes \id)(\id\otimes\sigma^{-1})(\sigma^{-1} 
\otimes\id)(m\otimes  \id^2)(\id\otimes  \sigma  \otimes  \id)(\phi 
\k\otimes \phi \k)\sigma_{-2} \\
=(\K^{-1}  \otimes  \k^{-1})(\ilG \otimes  \id)(\id\otimes  \sigma^{-1})
(\id\otimes m\otimes \id)(\sigma^{-1}\otimes \id^2 )(\phi\k\otimes\phi\k)
\sigma_{-2}\\
=(\K^{-1}\otimes\k^{-1})(\ilG \otimes m)(\id\otimes\sigma^{-1}\otimes 
\id)(\k\otimes \k\otimes \k\otimes \k)(\phi \otimes \phi )\sigma_{-2}\\
=(\irG\otimes m)[(\sigma_{-2})^{-1}\otimes (\sigma_{-2})^{-1}]
(\id\otimes \sigma^{-1}\otimes \id)(\phi \otimes \phi )\sigma_{-2}\\
=(\irG\otimes m)[(\sigma_{-2})^{-1}\otimes  \id^2]
(\id\otimes\phi\otimes\id)[\id\otimes (\sigma_{-1})^{-1}]
(\phi \otimes \id)\sigma_{-2}\\
=(\irG\otimes m)[(\sigma_{-2})^{-1}\otimes \id ](\id\otimes  \phi  \otimes 
\id)(\sigma_{-1}\otimes\id)(\id\otimes\phi)=(\irG\otimes 
m)(\id\otimes   \sigma\otimes\id)(\phi\otimes\phi).
\end{multline*}
Consequently $\Gamma$ is right-covariant with $\xi=\rac$ and
\eqref{59} holds.

Similarly, if $\Gamma$ is $\k$-covariant and  right-covariant
then a map $\xi\colon\Gamma\rightarrow\cal{A} \otimes \Gamma$ given by
$$
\xi =(\k^{-1}  \otimes \K^{-1})(\rtw{\sigma} )^{-1}\rac \K
$$
satisfies
\begin{equation*}
\begin{split}
 \xi\ilG&=(\k^{-1}  \otimes    \K^{-1})(\rtw{\sigma})^{-1}\rac
\irG(\k\otimes \k)\sigma_{-2}\\
&=(\k^{-1}  \otimes \K^{-1})(\rtw{\sigma} )^{-1}(\irG\otimes m)
(\id\otimes \sigma \otimes\id)(\phi \k\otimes \phi \k)\sigma_{-2}\\
&=(\k^{-1}  \otimes   \K^{-1})(m\otimes\irG)(\id\otimes\sigma^{-1}
\otimes\id)(\k\otimes \k\otimes \k\otimes \k)(\phi \otimes \phi )\sigma_{-2}\\
&=(m\otimes \ilG)[(\sigma_{-2})^{-1}\otimes(\sigma_{-2})^{-1}](\id\otimes 
\sigma^{-1}\otimes\id)(\phi \otimes     \phi )\sigma_{-2}\\&=(m\otimes 
\ilG)(\id\otimes  \sigma \otimes \id)(\phi \otimes \phi ).
\end{split}
\end{equation*}
This implies that $\Gamma$ is also left-covariant with $\xi=\lac$,
and equality~\eqref{510} holds. 

Finally,  let  us  assume  that $\Gamma$ is  bicovariant and
consider a map $\K\colon\Gamma\rightarrow\Gamma$ defined by
diagram \eqref{511}. Then a straightforward computation shows that
equality \eqref{51} holds, and that $\K$ is bijective.
In other words, $\Gamma$ is $\k$-covariant and
\eqref{511} holds by construction. 
\end{pf}

   Let us assume that $\Gamma$ is bicovariant. It turns out that
quadruplets $(\Gamma_{inv},\pi,{\circ},\cal{R})$ and
$(\Ginvr,\pir,\cirr,\cal{K})$ corresponding to the
left/right-covariant structure  on $\Gamma$  are naturally
related via the antipodal maps. According to \eqref{59}--\eqref{510}
\begin{equation}\label{512}
\K(\Gamma_{\inv})= \Ginvr\qquad \K(\Ginvr)=\Gamma_{\inv}.
\end{equation}

\begin{pro}\label{pro:53} The following identities hold
\begin{gather}
\K\pi =\pir  \k_0\qquad \quad\K\pir  =\pi \k_0 \label{513}\\
\K{\circ}(\pi \otimes \id)=\cirr(\k_0 \otimes \pir\k_0 )\tau\label{514a}\\ 
\K{\cirr}(\id\otimes \pir  )={\circ}(\pi \k_0 \otimes \k_0)\tau\label{514b}\\
\k_0(\cal{R})=\cal{K}\qquad \quad\k_0 (\cal{K})=\cal{R}\label{515}
\end{gather}
where $\k_0 =(\e\otimes \k)\sigma \phi =(\k\otimes \e)\sigma \phi$ is the
antipode associated to $\cal{A}_0$.
\end{pro}

\begin{pf} A direct computation gives
\begin{multline*}
\K\pi=\K\ilG(\k\otimes   \id)\phi= \irG\sigma_{-2}(\k^2\otimes\k)\phi 
= \irG(\id\otimes \k)\sigma_2(\k\otimes  \k)\phi= \mrG(d\otimes  \k)\sigma 
\tau^{-1}\phi\k\\
=\mrG\bigl[ \mrG(\pir   \otimes   \id)\phi   \otimes\k\bigr]
\sigma\tau^{-1}\phi\k= \mrG(\id\otimes m)(\pir  \otimes \id\otimes \k)
(\id\otimes \sigma  \tau^{-1})(\phi \otimes \id)\phi\\
= \mrG(\id\otimes m)(\pir  \otimes \id\otimes \k)(\id\otimes \sigma\tau^{-1} 
)(\id\otimes  \phi  )\phi= \mrG(\pir\otimes  1\e\k^{-1})\phi  \k=\pir 
(\k\otimes \e)\sigma \phi. 
\end{multline*}
Similarly,
\begin{multline*}
\K\pir=\K \irG(\id\otimes  \k)\phi   =\ilG \sigma_{-2}(\k\otimes\k^2)\phi 
=\ilG (\k\otimes \id)\sigma_2(\k\otimes \k)\sigma =\mlG (\k\otimes  d)\sigma 
\tau^{-1}   \phi \k\\ 
=\mlG \bigl[\k\otimes   \mlG (\id\otimes   \pi)\phi\bigr]\sigma\tau^{-1}\phi 
\k=\mlG (m\otimes \pi )(\k\otimes \id^2)(\sigma \tau^{-1}\otimes  \id)(\phi 
\otimes \id)\phi\\
=\mlG (1\e\k^{-1}  \otimes \pi )\phi =\pi (\e\otimes  \k)\sigma\phi
=\pi\k_0.
\end{multline*}
Relations \eqref{515} immediately  follow  from  \eqref{513},  definition  of 
spaces $\cal{R}$ and $\cal{K}$   and the fact that $\e\k_0=\e$. 

Let us check \eqref{514a}--\eqref{514b}. On the space 
$\ker(\e)\otimes \cal{A}$ the following equalities hold
$$
\K{\circ}(\pi \otimes\id)=\K\pi m\tau^{-1}\sigma=\pir\k_0 m_0
=\pir  m_0 (\k_0 \otimes\k_0 )\tau ={\cirr}(\k_0 \otimes \pir\k_0)\tau.
$$
Similarly, in the framework of the space $\cal{A} \otimes \ker(\e)$
we can write
$$
\K{\cirr}(\id\otimes \pir)=\K\pir  m_0 =\pi  \k_0 m_0 =\pi m_0
(\k_0 \otimes  \k_0 )\tau
={\circ}(\pi\k_0 \otimes\k_0)\tau .\qed
$$
\renewcommand{\qed}{}
\end{pf}

In terms of the bimodule  structuralizations  $\Gamma\leftrightarrow
\cal{A}  \otimes\Gamma_{\inv}$ and
$\Ginvr\otimes \cal{A}\leftrightarrow
\Gamma$ the operator $\K$ has a particularly simple form.

\begin{pro}\label{pro:54} The following identities hold
\begin{gather}
\K\mlG (\id\otimes \pi )= \mrG(\pir \k_0 \otimes \k)\tau\label{516a}\\
\K\mrG(\pir  \otimes \id)=\mlG (\k\otimes \pi\k_0 )\tau .\label{516b}
\end{gather}
\end{pro}

\begin{pf} We compute
\begin{multline*}
\K\mlG (\id\otimes \pi )= \mrG(\rtw{\sigma_{-2}} )(\k\otimes  \K\pi)
= \mrG( \rtw{\sigma_{-2}})(\k\otimes \pir  \k_0 )\\=
\mrG(\pir  \otimes  \id)\tau  (\k\otimes  \k_0 )=
\mrG(\pir  \k_0 \otimes \k)\tau .
\end{multline*}
Similarly,
\begin{multline*}
\K\mrG(\pir  \otimes \id)=(\mlG )\ltw{\sigma_{-2}}(\K\pir  \otimes\k)
=(\mlG)\ltw{\sigma_{-2}}(\pi \k_0 \otimes \k)\\
=\mlG(\id\otimes \pi)\tau(\k_0 \otimes \k)=\mlG(\k\otimes\pi\k_0)\tau .\qed
\end{multline*}
\renewcommand{\qed}{}
\end{pf}

\section{On $*$-Covariant Differential Structures}

Let us consider a quantum space $X$, represented by a unital
algebra $\cal{A}$ and assume that $X$ is $\cal{T}$-braided.
Let us also assume that $\cal{A}$ is equipped with a *-structure such that
\begin{equation}\label{62}
(*\otimes *)\alpha=\psi\alpha^{-1}\psi(*\otimes *),
\end{equation}
for each $\alpha\in\cal{T}$. Here, $\psi\colon\cal{A}
\otimes \cal{A}\rightarrow
\cal{A}\otimes\cal{A}$ is the standard transposition. It is 
easy to see that then \eqref{62} holds for every $\alpha\in\cal{T}^*$.
It is worth noticing that the operators
$$\cal{T}_c=\Bigl\{\psi\alpha^{-1}\psi\mid \alpha\in\cal{T}\Bigr\}$$
also form a braid system over $\cal{A}$. 

\begin{pro} Let us consider a *-covariant calculus $\Gamma$ over $X$. Then

\smallskip
\bla{i} If $\Gamma$ is left $\cal{T}$-covariant then it is also
left $\cal{T}_c$-covariant and
\begin{equation}\label{69}
\ltw{\alpha}(*\otimes *)=(*\otimes *)\ltw{(\psi\alpha^{-1}
\psi)},
\end{equation}
for each $\alpha\in\cal{T}$.

\smallskip
\bla{ii} Similarly, if $\Gamma$ is right $\cal{T}$-covariant then it is
right $\cal{T}_c$-covariant, with
\begin{equation}\label{610}
\rtw{\alpha}(*\otimes *)=(*\otimes *)\rtw{(\psi\alpha^{-1}\psi)},
\end{equation}
for each $\alpha\in\cal{T}$. \qed
\end{pro}

Let us now switch to multi-braided quantum groups $G$. Let us assume that
the *-structure on $\cal{A}$ satisfies
\begin{equation}\label{b32}
\phi*=(*\otimes*)\psi\sigma^{-1}\phi.
\end{equation}

\begin{defn} We shall say that the
antimultiplicative *-involution on $\cal{A}$ satisfying
the above equality is a *-structure on a braided quantum group $G$.
\end{defn}

This implies a number of further compatibility relations between $*$ and
maps appearing at the group level. At first, we have
\begin{equation}\label{b33}
\phi{*}\k=\bigl((*\k)\otimes(*\k)\bigr)\psi\phi. 
\end{equation}
The above equality implies
\begin{equation}\label{b34}
                   \e(a)^*=\e(\k(a)^*).
\end{equation}
Furthermore, as in the classical theory we have
\begin{equation}\label{b35}
                   \k^{-1}(a)=\k(a^*)^*
\end{equation}
for each $a\in\cal{A}$. Indeed, 
\begin{multline*}
\e(a)1=m\psi(*\k{*}\k\otimes\k)\psi\phi(a)=m\psi(\k\otimes
*\k{*}\k)\psi\phi(a)=\k[\k(a^{(1)})^*]^*\k(a^{(2)})\\
=\k(a^{(1)})\k[\k(a^{(2)})^*]^*,
\end{multline*}
and consequently
$$
a=\e(a^{(1)})a^{(2)}=\k[\k(a^{(1)})^*]^*\k(a^{(2)})a^{(3)}=
\k[\k(a^{(1)})^*]^*\e(a^{(2)})=\k[\k(a)^*]^*.
$$
Furthermore, let us examine interrelations between $*$, and braid
operators $\tau$ and $\sigma$.
\begin{pro} The following identities hold
\begin{align}
\sigma(*\otimes*)&=(*\otimes*)\psi\sigma^{-1}\psi\label{b36}\\
\tau(*\otimes*)&=(*\otimes*)\psi\tau^{-1}\psi.\label{b37}
\end{align}
\end{pro}

\begin{pf} Direct transformations give
\begin{multline*}    
\sigma(*\k\otimes*\k)=(*\otimes*)(m\psi\otimes m\psi)\bigl(\id\otimes\psi
\sigma^{-1}\phi m(\k\otimes\k)\psi\otimes\id\bigr)
(\psi\phi\otimes\psi\phi)\\
=(*\otimes *)\psi(m\otimes m)(\id\otimes\sigma^{-1}\phi m(\k\otimes\k)
\otimes\id)F\\ 
=(*\otimes*)\psi(m\otimes m)\bigl[\id\otimes\sigma^{-1}(m\otimes m)
(\id\otimes\sigma\otimes\id)(\phi\k\otimes\phi\k)\otimes\id\bigr]
F\\
=(*\otimes*)\psi(m\otimes m)(\id\otimes m\otimes m\otimes\id)
\bigl[\id\otimes(\id\otimes\sigma^{-1}
\otimes\id)(\sigma^{-1}\otimes\sigma^{-1})(\phi\k\otimes\phi\k)\otimes
\id\bigr]F\\
=(*\otimes*)\psi(m\otimes m)(m\otimes\sigma^{-1}\otimes m)(\id\otimes\k
\otimes\k\otimes\k\otimes\k\otimes\id)(\phi\otimes\id^2\otimes\phi)
F\\
=(*\otimes*)\psi(m\otimes m)(\id\otimes\sigma^{-1}\otimes\id)(1\e\otimes
\k\otimes\k\otimes 1\e)(\phi\otimes\phi)\psi=(*\otimes*)\psi\sigma^{-1}
\psi(\k\otimes\k),
\end{multline*}
where $F=(\phi\otimes\phi)\psi$. This
proves \eqref{b36}.  Furthermore, applying \eqref{b36} and the
definition of $\tau$ we obtain

\begin{multline*}
\tau(*\otimes*)=(*\k\otimes*\k\otimes *\e)(\id\otimes\psi\sigma\psi)
(\psi\phi\otimes\id)\psi\sigma^{-1}\psi(\k^{-1}\otimes\k^{-1})\\
=(*\k\otimes*\k\otimes *\e)(\id\otimes\psi\sigma\psi)(\psi\otimes\id)(
\id\otimes\psi)(\psi\otimes\id)(\id\otimes\phi)\sigma^{-1}\psi
(\k^{-1}\otimes\k^{-1})\\
=(*\k\otimes*\k\otimes *\e)(\id\otimes\psi\sigma)(\psi\otimes\id)
(\id\otimes\psi\phi)\sigma^{-1}\psi(\k^{-1}\otimes\k^{-1})\\
=(*\k\otimes*\k\otimes *\e)(\id\otimes\psi)(\psi\otimes\id)(\id\otimes
\psi)(\sigma\otimes\id)(\id\otimes\phi)\sigma^{-1}\psi
(\k^{-1}\otimes\k^{-1})\\
=(*\k\otimes*\k)(\e\otimes\psi\tau^{-1})(\phi\otimes\id)\psi
(\k^{-1}\otimes\k^{-1})=(*\otimes*)\psi\tau^{-1}\psi.
\end{multline*}
This completes the proof. 
\end{pf}

Condition \eqref{b32} says that the comultiplication $\phi$ is a 
hermitian map, if $\cal{A}\otimes\cal{A}$ is endowed with the *-structure
induced by $\sigma$ and $*\colon\cal{A}\rightarrow\cal{A}$. 

\begin{pro} Let us consider a *-covariant differential calculus $\Gamma$ over
$G$. 

\smallskip
\bla{i}  Let  us  assume  that  $\Gamma$  is,  in 
addition, left-covariant. Then
\begin{align}
\lac*&=\ltw{\sigma}\rtw{\psi}(*\otimes *)\lac\label{613}\\
*\pi&=-\pi{*}\k.\label{614}
\end{align}

\bla{ii}
Similarly, if $\Gamma$ is right-covariant then
\begin{align}
\rac *&=\rtw{\sigma} \ltw{\psi}(*\otimes *) \rac\label{615}\\
*\pir &=-\pir{*}\k.\label{616} 
\end{align}

\bla{iii} If $\Gamma$ is $\k$-covariant then
\begin{equation}\label{617}
\K(\vartheta^*)=\K^{-1}(\vartheta)^*
\end{equation}
for each $\vartheta\in\Gamma$. \qed
\end{pro}

For the end of this section,  let us 
characterize *-covariance of a left-covariant calculus in terms of
the  corresponding  right $\cal{A}_0$-ideal.
It turns out that the characterization is the same as in the standard
theory. 

\begin{pro}
Let $\Gamma$ be an arbitrary left-covariant calculus over $G$ and
$\cal{R}\subseteq\ker(\e)$ the associated right $\cal{A}_0$-ideal.
Then the calculus $\Gamma$ is *-covariant if and only if 
$\cal{R}$ is $*\k$-invariant.
\end{pro}

\begin{pf} If $\Gamma$ is *-covariant  then \eqref{614},
together with the definition of $\cal{R}$, implies that $\cal{R}$ is 
$*\k$-invariant.

Conversely, let us assume that $\cal{R}$ is  $*\k$-invariant.
Then  the  formula \eqref{614} consistently defines an antilinear
involution $*\colon\Gamma_{\inv}\rightarrow\Gamma_{\inv}$

   According to Proposition~\ref{pro:38} the formula 
$(a\vartheta)^* =\vartheta^*a^*$, 
where $a\in \cal{A}$ and $\vartheta\in\Gamma_{\inv}$ consistently 
defines an antilinear extension $*\colon\Gamma\rightarrow\Gamma$.
Applying the elementary transformations with $d$ and $\pi$ we obtain
\begin{equation*}
\begin{split}
(da)^*&=[a^{(1)}\pi(a^{(2)})]^*=-\pi(\k(a^{(2)})^*)a^{(1)*}\\
&=-\pi (\k(a^{(2)})^* )\k(\k(a^{(1)})^*)\\
&=-\k(\k(a^{(3)})^*)d[\k(a^{(2)})^*]\k(\k(a^{(1)})^*)\\
&=a^{(3)*}\k(a^{(2)})^*d(a^{(1)*})
=\e(a^{(2)})^* d(a^{(1)*})=d(a^*).
\end{split}
\end{equation*}
Consequently, $\Gamma$ is *-covariant.
\end{pf}

Let us observe that the above proof is
the same as in the standard theory \cite{W} (braidings are not included). 
A   similar   characterization   of   *-covariance    holds    for
right-covariant structures.

\appendix

\section{Right-Covariant Calculi}

 Let $\Gamma$  be  a  right-covariant  first  order  differential 
calculus over a braided quantum group $G$. 
The corresponding right action $\rac\colon\Gamma\rightarrow
\Gamma\otimes \cal{A}$
can be also characterized by 
\begin{equation}\label{A1}
\rac\ilG =(\ilG \otimes m)(\id\otimes\sigma\otimes\id)(\phi\otimes\phi).
\end{equation}
The right action map satisfies equalities
\begin{gather}
\rac\mrG=(\mrG\otimes m)(\id\otimes \sigma \otimes \id)( \rac 
\otimes \phi )\label{A2}\\
\rac d=(d\otimes \id)\phi\label{A3}\\
(\id\otimes \e)\rac=\id\label{A4}\\
( \rac \otimes \id) \rac =(\id\otimes \phi ) \rac.\label{A5}
\end{gather}

Every right-covariant calculus  is automatically
right $\cal{F}$-covariant.  In
particular,  the  flip-over  operator
$\rtw{\sigma}\colon\cal{A}   \otimes\Gamma\rightarrow
\Gamma\otimes\cal{A}$ is  determined  by the diagram
\begin{equation}\label{A6}
\begin{CD}
\cal{A}\otimes\cal{A}\otimes\Gamma\otimes\cal{A} @>{\mbox{$\k\otimes\rac\mlG
\otimes\k$}}>> \cal{A}\otimes\Gamma\otimes\cal{A}\otimes\cal{A}\\
@A{\mbox{$\phi\otimes\rac$}}AA  @VV{\mbox{$\mlG\otimes m$}}V\\
\cal{A}\otimes\Gamma
@>>{\mbox{$\quad\rtw{\sigma}\quad$}}> \Gamma\otimes\cal{A}
\end{CD}
\end{equation}
The operator $\rtw{\sigma}$  expresses the left multiplicativity  of
$\rac$, via the diagram
\begin{equation}\label{A7}
\begin{CD}
\cal{A}\otimes\cal{A}\otimes\Gamma\otimes\cal{A}
@>{\mbox{$\id\otimes\rtw{\sigma}
\otimes\id$}}>> \cal{A}\otimes\Gamma\otimes\cal{A}\otimes\cal{A}\\
@A{\mbox{$\phi\otimes\rac$}}AA  @VV{\mbox{$\mlG\otimes m$}}V\\
\cal{A}\otimes\Gamma
@>>{\mbox{$\quad\rac\mlG\quad$}}> \Gamma\otimes\cal{A}
\end{CD}
\end{equation}
The following twisting properties hold:
\begin{equation}\label{A8}
(\rac\otimes\id)\rtw{\sigma_{n+m}}=(\id\otimes \sigma_n)(\rtw{\sigma_m}\otimes
\id)(\id\otimes\rac). 
\end{equation}

Let $\Ginvr$ be  the  set  of  all  right-invariant  elements  of
$\Gamma$. Then the map $Q\colon\Gamma\rightarrow\Gamma$ defined by
\begin{equation}\label{A9}
Q= \mrG(\id\otimes \k) \rac
\end{equation}
projects $\Gamma$ onto $\Ginvr$. Moreover,
\begin{equation}\label{A10}
Q\irG=(Qd\otimes \e)\sigma^{-1}\tau .
\end{equation}
The composition
\begin{equation}\label{A11}
\pir  =Qd= \mrG(d\otimes \k)\phi\colon\cal{A}\rightarrow\Ginvr
\end{equation}
is surjective. All  flip-over  operators  $\rtw{\sigma_n}$
map $\cal{A}\otimes\Ginvr$ onto $\Ginvr\otimes\cal{A}$. Their
restrictions on this space are given by

\begin{equation}\label{A12}
\rtw{\sigma_n}(\id\otimes \pir)=(\pir  \otimes \id)\tau, 
\end{equation}
for each $n\in\Bbb{Z}$.

As a right $\cal{A}$-module,  the space $\Gamma$  is
naturally identificable with $\Ginvr\otimes \cal{A}$.
The isomorphism is induced by the multiplication map $\mrG$. 
Moreover,
\begin{equation}\label{A14}
[\mrG\restr\Ginvr\otimes \cal{A}]^{-1}=(Q\otimes \id) \rac.
\end{equation}

In terms of the structuralization $\Gamma\leftrightarrow
\Ginvr\otimes\cal{A}$, the following correspondences hold
\begin{gather}
\mrG \leftrightarrow\id\otimes m\qquad\rac
\leftrightarrow\id\otimes \phi\label{A16}\\
d\leftrightarrow (\pir  \otimes \id)\phi\label{A17}\\
\mlG\leftrightarrow(\cirr\otimes m)(\id\otimes\sinvr\otimes \id)(\phi\otimes 
\id^2).\label{A18}
\end{gather}

 Here, $\sinvr\colon\cal{A}   \otimes\Ginvr\rightarrow\Ginvr\otimes 
\cal{A}$  is the restriction of the operators $\rtw{\sigma_n}$,
and the map $\cirr\colon\cal{A}\otimes\Ginvr\rightarrow\Ginvr$ is given by
\begin{equation}\label{A19}
a\cirr\vartheta=Q(a\vartheta). 
\end{equation}
This map determines a left $\cal{A}_0$-module structure on 
$\Ginvr$. We have also
\begin{equation}\label{A20}
a\cirr\pir(b)=\pir  m_0(a\otimes b)-\pir(a)\e(b). 
\end{equation}

The space $\Gamma$ is also trivial as a left $\cal{A}$-module.
The corresponding isomorphism $\Gamma\leftrightarrow\cal{A}\otimes\Ginvr$
is induced by the product map, and explicitly
\begin{gather*}   
(\mlG\restr\cal{A} \otimes\Ginvr)=(\cirr\otimes\id)(\id\otimes 
\sinvr)(\phi \otimes \id)\colon\cal{A} \otimes 
   \Ginvr\rightarrow\Ginvr\otimes \cal{A}\\
(\mlG\restr\cal{A} \otimes\Ginvr)^{-1}=
(\k\otimes\cirr)(\phi \k^{-1}\otimes 
\id)(\sinvr)^{-1}\colon\Ginvr\otimes \cal{A}\rightarrow\cal{A} \otimes\Ginvr.
\end{gather*}

  In terms of the structuralization $\Gamma=\cal{A}\otimes\Ginvr$
the following correspondences hold:
\begin{gather}
  \mrG\leftrightarrow
[m\otimes\cirr(\k^{-1}  \otimes  \id)](\id\otimes  \sigma^{-1}\phi
\otimes \id)(\id\otimes (\sinvr)^{-1})\label{A21}\\
\mlG\leftrightarrow m\otimes \id\label{A22}\\
\rac\leftrightarrow(\id\otimes    \sinvr)(\phi \otimes \id)\label{A23} \\
-d\leftrightarrow (\k\otimes \pir)\phi \k^{-1}  =(\id\otimes \pir  \k^{-1}
)\sigma^{-1}\phi .\label{A24}
\end{gather}

The structure of every right-covariant calculus $\Gamma$
is completely determined by the
space $\cal{K}=\ker(\pir)\cap\ker(\e)$. This space is a left
$\cal{A}_0$-ideal satisfying
\begin{equation}\label{A25}
\tau(\cal{A}\otimes\cal{K})=\cal{K}\otimes\cal{A} .
\end{equation}

Conversely, let $\cal{K}\subseteq\ker(\e)$ be a left $\cal{A}_0$-ideal
such that equality \eqref{A25} holds.
The space $\Ginvr$ and maps $\pir$ and $\cirr$ can be
recovered as
\begin{gather}
\Ginvr=\ker(\e)/\cal{K}\qquad\pir(a)=[a-\e(a)]_{\cal{K}}\\ 
\cirr(a\otimes [b]_{\cal{K}})=[m\tau^{-1}\sigma (a\otimes b)]_{\cal{K}}.
\end{gather}
The whole right-covariant calculus $\Gamma$ is then constructed  with  the 
help of the above established correspondences.

\section{Elementary Properties of The Adjoint Action}

By definition, the adjoint action of $G$ onto itself is a 
linear map $\ad\colon\cal{A}\rightarrow\cal{A}\otimes\cal{A}$
defined by 
\begin{equation}
\ad=(\id\otimes    m)(\id\otimes\k\otimes\id)(\tau     \otimes 
\id)(\id\otimes \phi )\phi.\label{B1}
\end{equation}
\begin{lem}\label{lem:B1} The following identities hold
\begin{align}
(\id\otimes\e)\ad&=\id\label{B3}\\ 
(\id\otimes \phi)\ad&=(\ad\otimes \id)\ad.\label{B4}
\end{align}
In other words, $\ad$ is a counital and coassociative map. 
\end{lem}

\begin{pf} We compute
\begin{multline*}
(\id\otimes \e)\ad=(\id\otimes \e\otimes \e)(\id\otimes 
\sigma^{-1}\tau )(\tau \otimes \id)(\k\otimes \phi )\phi\\ 
=(\id\otimes   \e\otimes   \e)(\tau 
\otimes\id)(\id\otimes\phi)\sigma^{-1}\tau (\k\otimes \id)\phi 
=(\id\otimes\e\k)\sigma\phi=\k^{-1}(\id\otimes \e)\phi\k=\id.
\end{multline*}
Computation of the left-hand side of \eqref{B4} gives
\begin{multline*}
(\id\otimes m\otimes  m)(\id^2 \otimes  \sigma  \otimes  \id)(\id\otimes 
\phi \otimes \phi )(\tau \otimes \id)(\k\otimes \phi )\phi=\\
=(\id\otimes m\otimes m)(\tau \otimes \sigma \otimes  \id)(\id\otimes 
\tau \otimes \id^2 )(\phi  \otimes  \id\otimes  \phi  )(\k\otimes  \phi 
)\phi\\
=(\id\otimes  m\otimes  m)(\tau  \otimes  \id^3 )(\k\otimes   \id\otimes 
\sigma \otimes \id)(\id\otimes \tau \otimes \id^2 )(\id\otimes  \k\otimes 
\id^3 )(\sigma \otimes \phi \otimes \id)(\phi \otimes \phi )\phi\\
=(\id\otimes   m\otimes   m)(\tau   \otimes   \id\otimes    \k\otimes 
\id)(\k\otimes \id\otimes \tau \sigma^{-1}\tau  \otimes  \id)(\id\otimes 
\tau \otimes \id^2 )(\sigma \otimes  \phi  \otimes  \id)(\phi  \otimes 
\phi )\phi\\
=(\id\otimes   m\otimes   m)(\tau   \otimes   \id\otimes    \k\otimes 
\id)(\k\otimes \id^4 )(\phi \otimes \id^3 )(\id\otimes \tau  \sigma^{-1}\tau 
\otimes \id)(\sigma \otimes \id^2 )(\phi \otimes \phi )\phi\\
=(\id\otimes   m\otimes   m)(\tau   \otimes   \id\otimes    \k\otimes 
\id)(\k\otimes \id^4)(\phi \otimes \id^3)(\phi \otimes \id^2)(\tau \otimes 
\id)(\id\otimes \phi )\phi\\
=(\ad\otimes m)(\id\otimes \k\otimes \id)(\tau  \otimes  \id)(\id\otimes 
\phi )\phi =(\ad\otimes \id)\ad,
\end{multline*}
which completes the proof.
\end{pf}

Further useful identities are
\begin{lem} We have
\begin{align*}
(\e\otimes \id)\ad&=1\e\\
\bigl[\id\otimes  m(\id\otimes  \k)\otimes   \id\bigr]
(\ad\otimes\phi)\phi
&=(\id\otimes \k\otimes \id)(\tau \otimes \id)(\id\otimes \phi  )\phi. 
\end{align*}
\end{lem}

\begin{pf} Let us check the second identity. A direct computation gives
\begin{multline*}
\bigl[\id\otimes m(\id\otimes \k)\otimes \id\bigr](\ad\otimes\phi)\phi(a)=\\= 
\bigr(\id\otimes m(\id\otimes 
\k)\otimes\id\bigr)(\id\otimes m\otimes 
\id^2)(\k\otimes    \id^4)(\tau \otimes    \id^3)(a^{(1)}\otimes a^{(2)}
\otimes a^{(3)}\otimes a^{(4)}\otimes a^{(5)}) \\
=(\id\otimes  m\otimes  \id)\bigl(\tau   \otimes   m(\id\otimes   \k)\otimes 
\id\bigr)(\k(a^{(1)})\otimes a^{(2)}\otimes a^{(3)}\otimes a^{(4)}\otimes
a^{(5)})\\
=(\id\otimes m\otimes \id)(\tau\otimes    \id^2)(\k(a^{(1)})\otimes 
a^{(2)}\otimes 1\e(a^{(3)})\otimes a^{(4)})\\
=(\id\otimes  m\otimes  \id)(\id\otimes  \k\otimes   \id)(\tau   \otimes 
\id)(\id\otimes \phi )\phi (a).\qed
\end{multline*}
\renewcommand{\qed}{}
\end{pf}
Finally, let us study the twisting properties of the adjoint action.

\begin{lem} The following identities hold
\begin{gather}
(\id\otimes\ad)\sigma_m=(\sigma_m\otimes\id)(\id\otimes\sigma_n)
(\ad\otimes\id)\label{B7}\\
(\ad\otimes\id)\sigma_n=(\id\otimes\sigma_m)(\sigma_n\otimes\id)
(\id\otimes\ad).\label{B8}
\end{gather}
\end{lem}

\begin{pf} We compute
\begin{multline*}
(\sigma_m\otimes \id)(\id\otimes \sigma_n)(\ad\otimes \id)=\\
=(\sigma_m   \otimes  \id)(\id\otimes  \sigma_n)(\id\otimes   m\otimes 
\id)(\tau \otimes \id^2 )(\k\otimes \id^3)\bigl((\phi \otimes\id)\phi\otimes 
\id\bigr)\\
=(\sigma_m\otimes m)(\id\otimes  \sigma_n\otimes  \id)(\tau  \otimes 
\sigma_n  )(\k\otimes \id^3)\bigl((\phi \otimes \id)\phi \otimes \id\bigr)\\ 
=(\id^2 \otimes m)(\id\otimes\tau \otimes\id)(\id\otimes\k\otimes\id^2)
(\sigma_{-n}   \otimes\id^2)(\id\otimes \sigma_m\otimes \id)(\id^2\otimes
\sigma_n)\bigl((\phi \otimes\id)\phi \otimes \id\bigr)\\
=(\id^2\otimes m)(\id\otimes 
\tau \otimes \id)(\id\otimes \k\otimes \id^2)\bigl(\id\otimes 
(\phi \otimes \id)\phi\bigr)\sigma_m=(\id\otimes \ad)\sigma_m.
\end{multline*}
Furthermore, we have
\begin{multline*}
(\id\otimes\sigma_m)(\sigma_n\otimes\id)(\id\otimes \ad)=
(\id\otimes\sigma_m)(\sigma_n\otimes     m)(\id\otimes
\tau \otimes \id)(\id\otimes \k\otimes \id^2)\bigl(\id\otimes(\phi\otimes 
\id)\phi\bigr)\\
=(\id\otimes m\otimes\id)(\id^2\otimes\sigma_m)(\id\otimes\sigma_m\otimes 
\id)(\sigma_n\otimes 
\id^2)(\id\otimes  \tau \otimes  \id)(\id\otimes  \k\otimes   \id^2)
\bigl(\id\otimes (\phi \otimes \id)\phi\bigr)\\ 
=(\id\otimes m\otimes\id)(\tau \otimes\sigma_m)(\id\otimes\sigma_n\otimes
\id)(\sigma_m\otimes\id^2)(\id\otimes \k\otimes \id^2)\bigl(\id\otimes
(\phi \otimes \id)\phi\bigr)\\
=(\id\otimes m\otimes \id)(\tau \otimes \id^2)(\k\otimes \id^3)
(\id^2 \otimes\sigma_m)(\id\otimes\sigma_n\otimes   \id)(\sigma_{-m}\otimes 
\id^2)\bigl(\id\otimes (\phi \otimes \id)\phi\bigr)\\
=(\id\otimes  m\otimes  \id)(\tau \otimes  \id^2)(\k\otimes 
\id^3 )\bigl((\phi\otimes\id)\phi\otimes\id\bigr)\sigma_n=(\ad\otimes 
\id)\sigma_n.\qed
\end{multline*}
\renewcommand{\qed}{}
\end{pf}

\end{document}